\documentstyle{amsart}
\makeatletter
\newskip\dgARROWLENGTH  \dgARROWLENGTH=2.5em\relax
\newskip\dgHORIZPAD     \dgHORIZPAD=1em\relax
\newskip\dgVERTPAD      \dgVERTPAD=2ex\relax
\newskip\dgLABELOFFSET  \dgLABELOFFSET=.7ex\relax
\newcount\dgARROWPARTS  \dgARROWPARTS=4\relax
\newcommand{\dgeverynode}{\displaystyle}
\newcommand{\dgeverylabel}{\scriptstyle}
\newskip\dgDOTSPACING   \dgDOTSPACING=0.35em
\newskip\dgDOTSIZE      \dgDOTSIZE=1.5\fontdimen8\tenln
\newcount\dgMAXSQUARE   \dgMAXSQUARE=4\relax
\newcount\dgMINDBLSQ    \dgMINDBLSQ=10\relax
\newskip\dgCOLUMNWIDTH  \dgCOLUMNWIDTH=2em\relax
\chardef\f@ur=4
\@namedef{dgo@}{\let\dg@VECTOR=\vector
   \dg@LBLPOS=\dgARROWPARTS \divide\dg@LBLPOS\tw@}%
\@namedef{dgo@-}{\let\dg@VECTOR=\line}%
\@namedef{dgo@!}{\let\dg@VECTOR=\dg@novector}%
\@namedef{dgo@1}{\dg@LBLPOS=\@ne\relax}%
\@namedef{dgo@2}{\dg@LBLPOS=\tw@\relax}%
\@namedef{dgo@3}{\dg@LBLPOS=\thr@@\relax}%
\@namedef{dgo@4}{\dg@LBLPOS=\f@ur\relax}%
\@namedef{dgo@5}{\dg@LBLPOS=5\relax}%
\@namedef{dgo@6}{\dg@LBLPOS=6\relax}%
\@namedef{dgo@7}{\dg@LBLPOS=7\relax}%
\@namedef{dgo@8}{\dg@LBLPOS=8\relax}%
\@namedef{dgo@9}{\dg@LBLPOS=9\relax}%
\def\dgt@e{\dg@DX=\@ne \dg@DY=\z@ \dg@SIZE=\@ne}%
\def\dgt@w{\dg@DX=\m@ne \dg@DY=\z@ \dg@SIZE=\@ne}%
\def\dgt@n{\dg@DX=\z@ \dg@DY=\@ne \dg@SIZE=\@ne}%
\def\dgt@s{\dg@DX=\z@ \dg@DY=\m@ne \dg@SIZE=\@ne}%
\def\dgt@ne{\dg@DX=\@ne \dg@DY=\@ne \dg@SIZE=\@ne}%
\def\dgt@se{\dg@DX=\@ne \dg@DY=\m@ne \dg@SIZE=\@ne}%
\def\dgt@nw{\dg@DX=\m@ne \dg@DY=\@ne \dg@SIZE=\@ne}%
\def\dgt@sw{\dg@DX=\m@ne \dg@DY=\m@ne \dg@SIZE=\@ne}%
\def\dgt@nne{\dg@DX=\@ne \dg@DY=\tw@ \dg@SIZE=\@ne}%
\def\dgt@nnw{\dg@DX=\m@ne \dg@DY=\tw@ \dg@SIZE=\@ne}%
\def\dgt@sse{\dg@DX=\@ne \dg@DY=-\tw@ \dg@SIZE=\@ne}%
\def\dgt@ssw{\dg@DX=\m@ne \dg@DY=-\tw@ \dg@SIZE=\@ne}%
\def\dgt@ene{\dg@DX=\tw@ \dg@DY=\@ne \dg@SIZE=\tw@}%
\def\dgt@ese{\dg@DX=\tw@ \dg@DY=\m@ne \dg@SIZE=\tw@}%
\def\dgt@wnw{\dg@DX=-\tw@ \dg@DY=\@ne \dg@SIZE=\tw@}%
\def\dgt@wsw{\dg@DX=-\tw@ \dg@DY=\m@ne \dg@SIZE=\tw@}%
\def\dggeometry{
   \dg@ZTEMP=\dg@XGRID \multiply\dg@ZTEMP\tw@
   \ifnum\dg@YGRID=\z@ \dg@ZTEMP=\tw@
   \else \divide\dg@ZTEMP\dg@YGRID \fi
   \ifnum\dg@ZTEMP>\f@ur \dg@ZTEMP=\f@ur \fi
   \ifnum\dg@ZTEMP<\@ne \dg@ZTEMP=\@ne \fi
   \unitlength=2sp\relax
   \ifnum\dg@ZTEMP<\tw@
      \advance\dg@ZTEMP\@ne
      \multiply\unitlength\dg@YGRID
   \else
      \multiply\unitlength\dg@XGRID \divide\unitlength\dg@ZTEMP
   \fi
   \dg@XGRID=\dg@ZTEMP \dg@YGRID=\tw@
   \dg@rmcommondiv\tw@\dg@XGRID\dg@YGRID
   \divide\unitlength\dg@YGRID \divide\unitlength\@m\relax}
\@namedef{dgo@..}{\let\dg@VECTOR=\dg@dotvector}%

\def\dg@dotvector(#1,#2)#3{%
   \begingroup
   \dg@XTEMP=#1\relax \dg@YTEMP=#2\relax
   \let\dg@NDOTS=\dg@XEND \let\dg@DOTDIAM=\dg@WEND
   \dg@NDOTS=\unitlength \multiply\dg@NDOTS #3\relax
   \dg@ZTEMP=\dg@YTEMP \dg@changesign\dg@YTEMP\dg@ZTEMP
   \ifnum\dg@XTEMP>\z@
      \ifnum\dg@YTEMP>\dg@XTEMP
         \multiply\dg@NDOTS\dg@YTEMP \divide\dg@NDOTS\dg@XTEMP \fi
   \else\ifnum\dg@XTEMP<\z@
      \ifnum\dg@YTEMP>-\dg@XTEMP
         \multiply\dg@NDOTS\dg@YTEMP \divide\dg@NDOTS-\dg@XTEMP \fi
   \fi\fi
   \dg@YTEMP=\dg@ZTEMP
   \divide\dg@NDOTS\dgDOTSPACING
   \ifnum\dg@NDOTS>\z@\else \dg@NDOTS=\@ne \fi
   \dg@ZTEMP=\unitlength \multiply\dg@ZTEMP #3\relax
   \divide\dg@ZTEMP\dg@NDOTS
   \ifnum\dg@XTEMP=\z@
      \dg@changesign\dg@ZTEMP\dg@YTEMP \dg@YTEMP=\dg@ZTEMP
   \else
      \dg@changesign\dg@ZTEMP\dg@XTEMP
      \multiply\dg@YTEMP\dg@ZTEMP \divide\dg@YTEMP\dg@XTEMP
      \dg@XTEMP=\dg@ZTEMP
   \fi
   \divide\dg@XTEMP\unitlength \divide\dg@YTEMP\unitlength
   \begin{picture}(0,0)
      \dg@DOTDIAM=\dgDOTSIZE \divide\dg@DOTDIAM\unitlength
      \multiput(0,0)(\dg@XTEMP,\dg@YTEMP){\dg@NDOTS}{%
         \circle*{\dg@DOTDIAM}}%
      \multiply\dg@XTEMP\dg@NDOTS \multiply\dg@YTEMP\dg@NDOTS
      \put(\dg@XTEMP,\dg@YTEMP){\vector(#1,#2){0}}%
   \end{picture}%
   \endgroup}%
\ifx\@auxdefsloaded\relax
\else\let\@auxdefsloaded\relax

\def\newenvironment{%
   \@ifnextchar *{\@@newenv{\global\@ignoretrue}}{\@@newenv{}*}}

\def\@@newenv#1*#2{%
   \@ifnextchar [{\@newenv{#1}{#2}}{\@newenv{#1}{#2}[0]}}

\long\def\@newenv#1#2[#3]#4#5{%
   \expandafter\newcommand\csname#2\endcsname[#3]{#4}%
   \expandafter\long\expandafter\def\csname end#2\endcsname{#5#1}}

\def\renewenvironment{%
   \@ifnextchar *{%
      \@@renewenv{\global\@ignoretrue}}{\@@renewenv{}*}}

\def\@@renewenv#1*#2{%
   \@ifnextchar [{\@renewenv{#1}{#2}}{\@renewenv{#1}{#2}[0]}}

\long\def\@renewenv#1#2[#3]#4#5{%
   \expandafter\renewcommand\csname#2\endcsname[#3]{#4}%
   \expandafter\long\expandafter\def\csname end#2\endcsname{#5#1}}

\def\newoptcommand#1#2{%
   \@ifnextchar [{\@optargdef#1#2}{\@optargdef#1#2[1]}}

\def\renewoptcommand#1#2{%
   \edef\@tempa{\expandafter\@cdr\string#1\@nil}%
   \@ifundefined{\@tempa}{%
      \@latexerr{\string#1\space undefined}\@ehc}{}%
   \@ifnextchar [{\@reoptargdef#1#2}{\@reoptargdef#1#2[1]}}

\long\def\@optargdef#1#2[#3]#4{%
   \@ifdefinable #1{\@reoptargdef#1#2[#3]{#4}}}

\long\def\@reoptargdef#1#2[#3]#4{%
   \@tempcnta#3\relax \@tempcntb \@ne
   \let#1\relax \let\@tempa\relax
   \edef\@tempb{\long\def\csname\string#1\endcsname
      [\@tempa\the\@tempcntb]}%
   \advance\@tempcntb \@ne \advance\@tempcnta \m@ne
   \@whilenum\@tempcnta>0\do{%
      \edef\@tempb{\@tempb\@tempa\the\@tempcntb}%
      \advance\@tempcntb \@ne \advance\@tempcnta \m@ne}%
   \let\@tempa=##\@tempb{#4}%
   \def#1{\@ifnextchar [{\csname\string#1\endcsname}{%
      \csname\string#1\endcsname[#2]}}}

\def\newoptenvironment{%
   \@ifnextchar *{\@@newoptenv{\global\@ignoretrue}}{%
      \@@newoptenv{}*}}

\def\@@newoptenv#1*#2#3{%
   \@ifnextchar [{\@newoptenv{#1}{#2}{#3}}{%
      \@newoptenv{#1}{#2}{#3}[0]}}

\long\def\@newoptenv#1#2#3[#4]#5#6{%
   \expandafter\newoptcommand\csname#2\endcsname{#3}[#4]{#5}%
   \expandafter\long\expandafter\def\csname end#2\endcsname{#6#1}}

\def\renewoptenvironment{%
   \@ifnextchar *{\@@renewoptenv{\global\@ignoretrue}}{%
      \@@renewoptenv{}*}}

\def\@@renewoptenv#1*#2#3{%
   \@ifnextchar [{\@renewoptenv{#1}{#2}{#3}}{%
      \@renewoptenv{#1}{#2}{#3}[0]}}

\long\def\@renewoptenv#1#2#3[#4]#5#6{%
   \expandafter\renewoptcommand\csname#2\endcsname{#3}[#4]{#5}%
   \expandafter\long\expandafter\def\csname end#2\endcsname{#6#1}}

\newcounter{keepoptional}
\newcounter{optnestctr}

\newoptenvironment*{optional}{1}[1]{%
   \ifnum#1>\value{keepoptional}\relax
      \setcounter{optnestctr}{0}\@powerup\expandafter\@endopt\fi
   \ignorespaces}{}

\def\@powerup{\catcode`\{=12 \catcode`\}=12 \catcode`\\=12 \relax}
\def\@powerdown{\catcode`\{=1 \catcode`\}=2 \catcode`\\=0 \relax}
\begingroup \catcode`|=0 \catcode`[=1 \catcode`]=2 \@powerup
   |long|gdef|@endopt#1\end{optional}[%
      |@beginopt#1\begin{optional}*]%
   |long|gdef|@beginopt#1\begin{optional}[%
      |@ifstar[%
         |ifnum|value[optnestctr]>0|relax
            |addtocounter[optnestctr][-1]|let|@zzzap=|@endopt
         |else
            |@powerdown|end[optional]|let|@zzzap=|relax|fi
	 |@zzzap]%
      [
         |addtocounter[optnestctr][1]|@beginopt]]%
|endgroup

\fi

\newcount\dg@HORIZ      \newcount\dg@VERT
\newcount\dg@XLPAD      \newcount\dg@YBPAD
\newcount\dg@XRPAD      \newcount\dg@YTPAD
\newcount\dg@X          \newcount\dg@Y
\newcount\dg@XGRID      \newcount\dg@YGRID
\newcount\dg@SIZE
\newcount\dg@USERSIZE
\newcount\dg@DX         \newcount\dg@DY
\newcount\dg@XLBL       \newcount\dg@YLBL
\newcount\dg@XOFFSET    \newcount\dg@YOFFSET
\newcount\dg@LBLOFF
\newcount\dg@XLBLOFF    \newcount\dg@YLBLOFF
\newcount\dg@LBLPOS
\newcount\dg@XTEMP      \newcount\dg@YTEMP
\newcount\dg@ZTEMP
\newcount\dg@XNODE      \newcount\dg@YNODE
\newcount\dg@XEND       \newcount\dg@YEND
\newcount\dg@WEND       \newcount\dg@HEND
\newcount\dg@COUNT
\newbox\dg@NODEBOX
\newoptenvironment*{diagram}{}[1]{%
   \global\advance\dg@COUNT\@ne \typeout{[diagram \the\dg@COUNT:}%
   \let\node=\dg@node \let\\=\dg@cr \let\arrow=\dg@arrow
   \def\dg@BIGNODE{#1}%
   \ifx\dg@BIGNODE\@empty\else
      \setbox\dg@NODEBOX=\hbox{$\dgeverynode{#1}$}%
      \dg@XTEMP=\wd\dg@NODEBOX
      \dg@YTEMP=\ht\dg@NODEBOX \advance\dg@YTEMP\dp\dg@NODEBOX
      \ifnum\dg@YTEMP=\z@ \dg@ZTEMP=\z@
      \else \dg@ZTEMP=\dg@XTEMP \divide\dg@ZTEMP\dg@YTEMP\fi
      \ifnum\dg@ZTEMP>\dgMAXSQUARE
         \ifnum\dg@ZTEMP<\dgMINDBLSQ
		 	\dg@XGRID=\thr@@ \dg@YGRID=\tw@
         \else \dg@XGRID=\tw@ \dg@YGRID=\@ne \fi
      \else \dg@XGRID=\@ne \dg@YGRID=\@ne \fi
      \advance\dg@XTEMP\dgHORIZPAD \advance\dg@YTEMP\dgVERTPAD
      \dg@XLPAD=\dg@XTEMP \divide\dg@XLPAD\tw@
      \dg@XRPAD=\dg@XTEMP \divide\dg@XRPAD\tw@
      \dg@YBPAD=\dg@YTEMP \divide\dg@YBPAD\tw@
      \dg@YTPAD=\dg@YTEMP \divide\dg@YTPAD\tw@
      \advance\dg@XTEMP\dgARROWLENGTH
      \divide\dg@XTEMP\dg@XGRID \divide\dg@XTEMP\@m
      \unitlength=1sp\relax \multiply\unitlength\dg@XTEMP
   \fi
   \typeout{saving...}%
   \dg@X=-\@ne \dg@Y=\z@ \dg@HORIZ=\z@\relax
   \let\dg@SLIST=\@empty
   \let\dg@NLIST=\@empty \let\dg@ALIST=\@empty
   \let\dg@PASS=\dg@savepass
}{
   \dg@VERT=-\dg@Y\relax
   \typeout{calculating...}%
   \ifx\dg@BIGNODE\@empty
      \dg@XGRID=\z@ \dg@YGRID=\z@
      \dg@XLPAD=\z@ \dg@XRPAD=\z@ \dg@YBPAD=\z@ \dg@YTPAD=\z@
      \let\dg@PASS=\dg@geompass
      \dg@NLIST \dg@ALIST
      \ifnum\dg@XGRID>\z@\else \dg@XGRID=\dgHORIZPAD \fi
      \ifnum\dg@YGRID>\z@\else \dg@YGRID=\dgVERTPAD \fi
      \dggeometry
      \ifdim\unitlength>\z@\else
         \unitlength=\dgARROWLENGTH \divide\unitlength\@m \fi
      \ifnum\dg@XGRID>\z@\else \dg@XGRID=\@ne \fi
      \ifnum\dg@YGRID>\z@\else \dg@YGRID=\@ne \fi
   \fi
   \dg@LBLOFF=\dgLABELOFFSET \divide\dg@LBLOFF\unitlength
   \multiply\dg@HORIZ\@m \multiply\dg@HORIZ\dg@XGRID
   \multiply\dg@VERT\@m \multiply\dg@VERT\dg@YGRID
   \divide\dg@XLPAD\unitlength \divide\dg@XRPAD\unitlength
   \divide\dg@YBPAD\unitlength \divide\dg@YTPAD\unitlength
   \typeout{drawing...}%
   \let\dg@PASS=\dg@drawpass
   \hspace*{\dg@XLPAD\unitlength}%
   \vcenter{%
      \hsize=0pt\relax\parindent=0pt\relax
      \parskip=0pt\relax\baselineskip=0pt\relax
      \vspace*{\dg@YTPAD\unitlength}%
      \begin{picture}(0,0)(0,0)\dg@NLIST\dg@ALIST\end{picture}%
      \raisebox{\z@}[\z@][\dg@VERT\unitlength]{}%
      \vspace*{\dg@YBPAD\unitlength}%
   }
   \hspace*{\dg@HORIZ\unitlength}%
   \hspace*{\dg@XRPAD\unitlength}%
   \typeout{done]}}%

\def\dg@savepass{s}
\def\dg@geompass{g}
\def\dg@drawpass{d}

\newoptcommand{\dg@node}{\@ne}[2]{%
   \ifx\dg@PASS\dg@savepass
      \dg@XTEMP=#1\relax
      \ifnum\dg@XTEMP<\@ne \dg@XTEMP=\@ne\fi
      \advance\dg@X\dg@XTEMP
      \ifnum\dg@HORIZ<\dg@X \dg@HORIZ=\dg@X \fi
      \setbox\dg@NODEBOX=\hbox{$\dgeverynode{#2}$}%
      \dg@XTEMP=\wd\dg@NODEBOX \advance\dg@XTEMP\dgHORIZPAD
      \dg@YTEMP=\ht\dg@NODEBOX \advance\dg@YTEMP\dp\dg@NODEBOX
      \advance\dg@YTEMP\dgVERTPAD
      \toks\z@=\expandafter{\dg@SLIST}%
      \edef\dg@SLIST{\the\toks\z@
         ,\noexpand\dg@XNODE=\number\dg@X\noexpand\relax
         \noexpand\dg@YNODE=\number\dg@Y\noexpand\relax
         \noexpand\dg@XTEMP=\number\dg@XTEMP\noexpand\relax
         \noexpand\dg@YTEMP=\number\dg@YTEMP\noexpand\relax}%
      \toks\z@=\expandafter{\dg@NLIST}%
      \toks\tw@={\dg@node{#2}}%
      \edef\dg@NLIST{\the\toks\z@
         \noexpand\dg@X=\number\dg@X\noexpand\relax
         \noexpand\dg@Y=\number\dg@Y\noexpand\relax
         \the\toks\tw@}%
   \else\ifx\dg@PASS\dg@geompass
      \ifnum\dg@X=\z@
         \dg@getnodesize
            {\dg@SLIST}{\dg@X}{\dg@Y}{\dg@WEND}{\dg@HEND}%
         \divide\dg@WEND\tw@
         \ifnum\dg@XLPAD<\dg@WEND \dg@XLPAD=\dg@WEND \fi\fi
      \ifnum\dg@X=\dg@HORIZ
         \dg@getnodesize
            {\dg@SLIST}{\dg@X}{\dg@Y}{\dg@WEND}{\dg@HEND}%
         \divide\dg@WEND\tw@
         \ifnum\dg@XRPAD<\dg@WEND \dg@XRPAD=\dg@WEND \fi\fi
      \ifnum\dg@Y=\z@
         \dg@getnodesize
            {\dg@SLIST}{\dg@X}{\dg@Y}{\dg@WEND}{\dg@HEND}%
         \divide\dg@HEND\tw@
         \ifnum\dg@YTPAD<\dg@HEND \dg@YTPAD=\dg@HEND \fi\fi
      \ifnum\dg@Y=-\dg@VERT\relax
         \dg@getnodesize
            {\dg@SLIST}{\dg@X}{\dg@Y}{\dg@WEND}{\dg@HEND}%
         \divide\dg@HEND\tw@
         \ifnum\dg@YBPAD<\dg@HEND \dg@YBPAD=\dg@HEND \fi\fi
   \else\ifx\dg@PASS\dg@drawpass
      \dg@XNODE=\dg@X \multiply\dg@XNODE\@m
      \multiply\dg@XNODE\dg@XGRID
      \dg@YNODE=\dg@Y \multiply\dg@YNODE\@m
      \multiply\dg@YNODE\dg@YGRID
      \setbox\dg@NODEBOX=\hbox{$\dgeverynode{#2}$}%
      \put(\dg@XNODE,\dg@YNODE){\dg@makebox{\box\dg@NODEBOX}}%
   \fi\fi\fi}%

\newoptcommand{\dg@cr}{\@ne}[1]{%
   \ifx\dg@PASS\dg@savepass
      \dg@YTEMP=#1\relax
      \ifnum\dg@YTEMP<\@ne \dg@YTEMP=\@ne \fi
      \advance\dg@Y -\dg@YTEMP\relax
      \dg@X=-\@ne\relax\fi}%
\newoptcommand{\dg@arrow}{\@ne}[2]{%
   \begingroup
   \dg@USERSIZE=#1\relax
   \ifnum\dg@USERSIZE<\@ne \dg@USERSIZE=\@ne \fi
   \dg@parse{#2}%
   \ifx\dg@PASS\dg@savepass
      \ifx\dg@label\dgl@b \let\dg@label=\dgl@t \fi
      \ifx\dg@label\dgl@r \let\dg@label=\dgl@l \fi
      \let\dg@process=\dg@save
   \else\ifx\dg@PASS\dg@geompass
      \let\dg@process=\dg@ignore
      \dg@geomcalc
   \else\ifx\dg@PASS\dg@drawpass
      \let\dg@process=\dg@draw
      \dg@drawcalc
   \fi\fi\fi
   \dg@label{\dg@process{#1}{#2}}}%

\newoptcommand{\arrow}{\@ne}[2]{%
   \dg@parse{#2}%
   \ifx\dg@label\dgl@b \let\dg@label=\dgl@t \fi
   \ifx\dg@label\dgl@r \let\dg@label=\dgl@l \fi
   \dg@label{\dg@textarrow{#1}{#2}}}%

\def\dg@textarrow#1#2#3#4{%
   \mathop{{\dgHORIZPAD=0pt\relax\dgVERTPAD=0pt\relax
      \begin{diagram}
         \node{}\arrow[#1]{#2}{#3}{#4}\node{}
      \end{diagram}}}}

\def\dg@parse#1{%
   \let\dg@label=\dgl@ \dgo@
   \let\dg@type=\@empty \let\dg@lbltype=\@empty
   \@for\dg@list:=#1\do{%
      \ifx\dg@type\@empty \let\dg@type=\dg@list
      \else\ifx\dg@lbltype\@empty \let\dg@lbltype=\dg@list
         \@ifundefined{dgo@\dg@list}{}{\@nameuse{dgo@\dg@list}}%
      \else
         \@ifundefined{dgo@\dg@list}{}{\@nameuse{dgo@\dg@list}}%
      \fi\fi}%
   \@ifundefined{dgt@\dg@type}{\dgt@e}{\@nameuse{dgt@\dg@type}}%
   \@ifundefined{dgl@\dg@lbltype}{}{%
      \dg@letname\dg@label{dgl@\dg@lbltype}}}

\def\dg@draw#1#2#3#4{%
   \put(\dg@X,\dg@Y){\dg@makebox{%
      \begin{picture}(0,0)%
         \thinlines
         \put(\dg@XOFFSET,\dg@YOFFSET){%
            \dg@VECTOR(\dg@DX,\dg@DY){\dg@SIZE}}%
         \put(\dg@XLBL,\dg@YLBL){\dg@makebox{%
            \begin{picture}(0,0)%
               \put(\dg@XLBLOFF,\dg@YLBLOFF){%
                  \dg@makebox[\dg@LBLONE]{$\dgeverylabel{#3}$}}%
               \put(-\dg@XLBLOFF,-\dg@YLBLOFF){%
                  \dg@makebox[\dg@LBLTWO]{$\dgeverylabel{#4}$}}%
            \end{picture}}}%
      \end{picture}}}%
   \endgroup}%

\def\dg@save#1#2#3#4{%
   \endgroup 
   \toks\z@=\expandafter{\dg@ALIST}%
   \toks\tw@={\dg@arrow[#1]{#2}{#3}{#4}}%
   \edef\dg@ALIST{\the\toks\z@%
      \noexpand\dg@X=\number\dg@X\noexpand\relax
      \noexpand\dg@Y=\number\dg@Y\noexpand\relax
      \the\toks\tw@}}%

\def\dg@ignore#1#2#3#4{\endgroup}

\def\dg@geomcalc{%
   \dg@XEND=\dg@SIZE \multiply\dg@XEND\dg@USERSIZE
   \ifnum\dg@DX=\z@
      \dg@YEND=\dg@XEND \dg@XEND=\z@
      \dg@changesign\dg@YEND\dg@DY
   \else
      \dg@changesign\dg@XEND\dg@DX \dg@YEND=\dg@XEND
      \multiply\dg@YEND\dg@DY \divide\dg@YEND\dg@DX
   \fi
   \advance\dg@XEND\dg@X \advance\dg@YEND\dg@Y
   \dg@getnodesize
      {\dg@SLIST}{\dg@XEND}{\dg@YEND}{\dg@WEND}{\dg@HEND}%
   \dg@XOFFSET=\dg@WEND \dg@YOFFSET=\dg@HEND
   \dg@getnodesize
      {\dg@SLIST}{\dg@X}{\dg@Y}{\dg@WEND}{\dg@HEND}%
   \advance\dg@XOFFSET\dg@WEND \divide\dg@XOFFSET\tw@
   \advance\dg@YOFFSET\dg@HEND \divide\dg@YOFFSET\tw@
   \dg@XTEMP=\dgARROWLENGTH \dg@YTEMP=\dgARROWLENGTH
   \ifnum\dg@DX>\z@
      \dg@ZTEMP=\dg@DX \multiply\dg@XTEMP\dg@DX
   \else \dg@ZTEMP=-\dg@DX \multiply\dg@XTEMP -\dg@DX \fi
   \ifnum\dg@DY>\z@
      \advance\dg@ZTEMP\dg@DY \multiply\dg@YTEMP\dg@DY
   \else \advance\dg@ZTEMP -\dg@DY \multiply\dg@YTEMP -\dg@DY\fi
   \ifnum\dg@ZTEMP=\z@\else
      \divide\dg@XTEMP\dg@ZTEMP \divide\dg@YTEMP\dg@ZTEMP
      \advance\dg@XOFFSET\dg@XTEMP \advance\dg@YOFFSET\dg@YTEMP
   \fi
   \divide\dg@XOFFSET\dg@SIZE \divide\dg@XOFFSET\dg@USERSIZE
   \divide\dg@YOFFSET\dg@SIZE \divide\dg@YOFFSET\dg@USERSIZE
   \ifnum\dg@DX=\z@ \dg@XOFFSET=\z@ \fi
   \ifnum\dg@DY=\z@ \dg@YOFFSET=\z@ \fi
   \ifnum\dg@XGRID<\dg@XOFFSET \global\dg@XGRID=\dg@XOFFSET\fi
   \ifnum\dg@YGRID<\dg@YOFFSET \global\dg@YGRID=\dg@YOFFSET\fi
   \relax}

\def\dg@drawcalc{%
   \dg@XEND=\dg@SIZE \multiply\dg@XEND\dg@USERSIZE
   \ifnum\dg@DX=\z@
      \dg@YEND=\dg@XEND \dg@XEND=\z@
      \dg@changesign\dg@YEND\dg@DY
   \else
      \dg@changesign\dg@XEND\dg@DX \dg@YEND=\dg@XEND
      \multiply\dg@YEND\dg@DY \divide\dg@YEND\dg@DX
   \fi
   \advance\dg@XEND\dg@X \advance\dg@YEND\dg@Y
   \dg@getnodesize
      {\dg@SLIST}{\dg@XEND}{\dg@YEND}{\dg@WEND}{\dg@HEND}%
   \divide\dg@WEND\unitlength \divide\dg@HEND\unitlength
   \multiply\dg@DX\dg@XGRID \multiply\dg@DY\dg@YGRID
   \dg@rmcommondiv\tw@\dg@DX\dg@DY
   \dg@rmcommondiv\tw@\dg@DX\dg@DY 
   \dg@rmcommondiv\thr@@\dg@DX\dg@DY
   \multiply\dg@SIZE\dg@USERSIZE \multiply\dg@SIZE\@m
   \ifnum\dg@DX=\z@
      \multiply\dg@SIZE\dg@YGRID
      \divide\dg@HEND\tw@ \advance\dg@SIZE -\dg@HEND
      \dg@getnodesize
         {\dg@SLIST}{\dg@X}{\dg@Y}{\dg@WEND}{\dg@YOFFSET}%
      \divide\dg@YOFFSET\unitlength \divide\dg@YOFFSET\tw@
      \advance\dg@SIZE -\dg@YOFFSET
      \dg@XOFFSET=\z@
      \def\dg@LBLONE{r}\def\dg@LBLTWO{l}%
      \dg@XLBL=\z@ \dg@YLBL=\dg@SIZE
      \multiply\dg@YLBL\dg@LBLPOS
      \divide\dg@YLBL\dgARROWPARTS\relax
      \advance\dg@YLBL\dg@YOFFSET
      \dg@changesign\dg@YLBL\dg@DY
      \dg@changesign\dg@YOFFSET\dg@DY
   \else
      \multiply\dg@SIZE\dg@XGRID
      \ifnum\dg@DY=\z@
         \divide\dg@WEND\tw@ \advance\dg@SIZE -\dg@WEND
         \dg@getnodesize
            {\dg@SLIST}{\dg@X}{\dg@Y}{\dg@XOFFSET}{\dg@HEND}%
         \divide\dg@XOFFSET\unitlength \divide\dg@XOFFSET\tw@
         \advance\dg@SIZE -\dg@XOFFSET
         \dg@YOFFSET=\z@
         \def\dg@LBLONE{b}\def\dg@LBLTWO{t}%
         \dg@YLBL=\z@ \dg@XLBL=\dg@SIZE
         \multiply\dg@XLBL\dg@LBLPOS
         \divide\dg@XLBL\dgARROWPARTS\relax
         \advance\dg@XLBL\dg@XOFFSET
         \dg@changesign\dg@XLBL\dg@DX
         \dg@changesign\dg@XOFFSET\dg@DX
      \else
         \divide\dg@WEND\tw@ \divide\dg@HEND\tw@
         \multiply\dg@HEND\dg@DX \divide\dg@HEND\dg@DY
         \ifnum\dg@HEND<\z@ \multiply\dg@HEND\m@ne \fi
         \ifnum\dg@WEND<\dg@HEND \advance\dg@SIZE -\dg@WEND
         \else \advance\dg@SIZE -\dg@HEND \fi
         \dg@getnodesize
            {\dg@SLIST}{\dg@X}{\dg@Y}{\dg@WEND}{\dg@HEND}%
         \divide\dg@WEND\unitlength \divide\dg@WEND\tw@
         \divide\dg@HEND\unitlength \divide\dg@HEND\tw@
         \multiply\dg@HEND\dg@DX \divide\dg@HEND\dg@DY
         \ifnum\dg@HEND<\z@ \multiply\dg@HEND\m@ne \fi
         \ifnum\dg@WEND<\dg@HEND \dg@XOFFSET=\dg@WEND
         \else \dg@XOFFSET=\dg@HEND \fi
         \advance\dg@SIZE -\dg@XOFFSET
         \dg@changesign\dg@XOFFSET\dg@DX
         \dg@YOFFSET=\dg@XOFFSET
         \multiply\dg@YOFFSET\dg@DY \divide\dg@YOFFSET\dg@DX
         \def\dg@LBLONE{br}\def\dg@LBLTWO{tl}%
         \ifnum\dg@DX<\z@ \ifnum\dg@DY>\z@
            \def\dg@LBLONE{bl}\def\dg@LBLTWO{tr}\fi\fi
         \ifnum\dg@DX>\z@ \ifnum\dg@DY<\z@
            \def\dg@LBLONE{bl}\def\dg@LBLTWO{tr}\fi\fi
         \dg@XLBL=\dg@SIZE
         \multiply\dg@XLBL\dg@LBLPOS
         \divide\dg@XLBL\dgARROWPARTS\relax
         \dg@changesign\dg@XLBL\dg@DX
         \dg@YLBL=\dg@XLBL
         \multiply\dg@YLBL\dg@DY \divide\dg@YLBL\dg@DX
         \advance\dg@XLBL\dg@XOFFSET
         \advance\dg@YLBL\dg@YOFFSET
      \fi
   \fi
   \dg@XLBLOFF=-\dg@DY \dg@changesign\dg@XLBLOFF\dg@DX
   \dg@YLBLOFF=\dg@DX \dg@changesign\dg@YLBLOFF\dg@DX
   \ifnum\dg@DX=\z@ \dg@XLBLOFF=\m@ne \fi
   \dg@XTEMP=\dg@DX \dg@changesign\dg@XTEMP\dg@DX
   \dg@YTEMP=\dg@DY \dg@changesign\dg@YTEMP\dg@DY
   \ifnum\dg@YTEMP>\dg@XTEMP \dg@XTEMP=\dg@YTEMP \fi
   \ifnum\dg@XTEMP=\z@ \dg@XTEMP=\@ne \fi
   \multiply\dg@XLBLOFF\dg@LBLOFF \divide\dg@XLBLOFF\dg@XTEMP
   \multiply\dg@YLBLOFF\dg@LBLOFF \divide\dg@YLBLOFF\dg@XTEMP
   \multiply\dg@X\@m \multiply\dg@X\dg@XGRID
   \multiply\dg@Y\@m \multiply\dg@Y\dg@YGRID
   \relax}%

\def\dg@rmcommondiv#1#2#3{%
   \dg@XTEMP=#2\relax
   \divide\dg@XTEMP #1\relax \multiply\dg@XTEMP #1\relax
   \dg@YTEMP=#3\relax
   \divide\dg@YTEMP #1\relax \multiply\dg@YTEMP #1\relax
   \ifnum\dg@XTEMP=#2\relax \ifnum\dg@YTEMP=#3\relax
      \divide#2#1\relax \divide#3#1\relax \fi\fi}%

\def\dg@changesign#1#2{%
   \ifnum #2<\z@ \multiply#1\m@ne
   \else\ifnum #2=\z@ #1=\z@ \fi\fi}%

\def\dg@getnodesize#1#2#3#4#5{%
   #4=\z@\relax #5=\z@\relax
   \expandafter\@for\expandafter\dg@trynode
   \expandafter:\expandafter=#1\do{%
      \dg@XNODE=\m@ne 
      \dg@trynode
      \ifnum #2=\dg@XNODE \ifnum #3=\dg@YNODE
         #4=\dg@XTEMP\relax #5=\dg@YTEMP\relax\fi\fi}}%

\long\def\dg@sinkbaseline#1{%
   \leavevmode\hbox{\vbox{%
      \lineskiplimit=16383pt\relax\lineskip=0pt\relax
      \baselineskip=-1000pt\relax
      \parindent=0pt\relax\parskip=0pt\relax
      \hbox{#1}\rule{0pt}{0pt}}}}%
\newoptcommand{\dg@makebox}{}[2]{\dg@sinkbaseline{%
   \expandafter\makebox\expandafter(\expandafter
      0\expandafter,\expandafter0\expandafter)\expandafter
      [#1]{#2}}}%

\def\dg@novector(#1,#2)#3{}%

\def\dg@letname#1#2{%
   \relax\expandafter
   \let\expandafter #1\csname #2\endcsname\relax}%

\def\dgl@#1{#1{}{}}%
\def\dgl@t#1#2{#1{#2}{}}%
\def\dgl@b#1#2{#1{}{#2}}%
\def\dgl@tb#1#2#3{#1{#2}{#3}}%
\def\dgl@l#1#2{#1{#2}{}}%
\def\dgl@r#1#2{#1{}{#2}}%
\def\dgl@lr#1#2#3{#1{#2}{#3}}%
\makeatother

\expandafter\ifx\csname >amspapermacs\endcsname\relax
\expandafter\def\csname >amspapermacs\endcsname{done}
\else\endinput\fi
\newtheorem{thm}{Theorem}[section]
\newtheorem{cor}[thm]{Corollary}
\newtheorem{lem}[thm]{Lemma}
\newtheorem{prop}[thm]{Proposition}

\theoremstyle{definition}
\newtheorem{definition}[thm]{Definition}

\theoremstyle{remark}
\newtheorem{remark}[thm]{Remark}
\newtheorem{ex}[thm]{Example}

\numberwithin{equation}{section}

\newcommand{\thmref}[1]{Theorem~\ref{#1}}
\newcommand{\secref}[1]{Section~\ref{#1}}
\newcommand{\lemref}[1]{Lemma~\ref{#1}}
\newcommand{\propref}[1]{Prop\-o\-si\-tion~\ref{#1}}
\newcommand{\corref}[1]{Cor\-ol\-lary~\ref{#1}}

\makeatletter
\newif\ifdraft
\draftfalse

\let\foob@r=\label
\let\f@@index\index
\def\t@index#1{\relax}
\def\label#1{\foob@r{#1}\ifdraft
\index{\string\labelentry{\@currenvir}{#1}}
\ifinner
\else\leavevmode\marginpar[\hfill
{\tiny\normalshape #1}]{{\tiny\normalshape
#1}}\fi\fi\ignorespaces}
\def\@lem@#1{Lemma~\ref{#1}}
\def\@equation@#1{Equation~(\ref{#1})}
\let\@align@=\@equation@
\let\@gather@=\@equation@
\def\@thm@#1{Theorem~\ref{#1}}
\def\@cor@#1{Corollary~\ref{#1}}
\def\@prop@#1{Proposition~\ref{#1}}
\def\@definition@#1{Definition~\ref{#1}}
\def\@remark@#1{Remark~\ref{#1}}
\def\@ex@#1{Example~\ref{#1}}
\def\dotfill{\leaders\hbox to0.5em{\hss.\hss}\hfill}
\newdimen\lcrskip
\def\labelentry#1#2#3{\vskip0pt plus 0.25pt\hbox
to\hsize{\hskip\lcrskip\@ifundefined{@#1@}{{\tt{\uppercase{#1}}}~\ref{#2}}
{\csname @#1@\endcsname{#2}}\dotfill{#2}\quad
(Page #3)\hskip\lcrskip}}
\def\labeltable{\ifdraft\clearpage\immediate\closeout\@indexfile
\section*{Label Cross References}
\input \jobname.idx
\else\relax\fi}

\def\gmeasure@#1{\let\f@@index\index\let\index\t@index\gwidth@
\z@\gmaxwidth@\z@\setbox@ne\vbox{\Let@
 \firstchoice@false\let\tag\gobble@tag
 \halign{\setboxz@h{$\m@th\displaystyle{##}$}\global\gwidth@\wdz@
\ifdim\gwidth@>\gmaxwidth@\global\gmaxwidth@\gwidth@\fi
 &\@gobble{##}\crcr#1\crcr}}\let\index\f@@index}
\def\measure@#1{\let\f@@index\index\let\index\t@index
\lwidth@\z@\rwidth@\z@\maxlwidth@\z@\maxrwidth@\z@
 \global\and@\z@
 \setbox@ne\vbox{%
   \everycr{\noalign{\global\tag@false\global\and@\z@}}\Let@
   \let\tag\gobble@tag
   \let\notag\@empty \let\nonumber\@empty
   \firstchoice@false
    \halign{\setboxz@h{$\m@th\displaystyle{\@lign##}$}%
     \global\lwidth@\wdz@
     \ifdim\lwidth@>\maxlwidth@\global\maxlwidth@\lwidth@\fi
     \global\advance\and@\@ne
     &\setboxz@h{$\m@th\displaystyle{{}\@lign##}$}%
     \global\rwidth@\wdz@
     \ifdim\rwidth@>\maxrwidth@\global\maxrwidth@\rwidth@\fi
     \global\advance\and@\@ne
     &\Tag@\@gobble{##}\crcr#1\crcr}}%
 \totwidth@\maxlwidth@\advance\totwidth@\maxrwidth@
\let\index\f@@index}
\newbox\di@b@x

\newbox\di@b@x
\newenvironment{ediagram*}{\setbox\di@b@x=\hbox\bgroup$
\begin{diagram}}{\end{diagram}$\egroup%
\begin{equation*}\copy\di@b@x\end{equation*}\global\@ignoretrue}
\makeatother

\newcount\hours
\newcount\minutes       
\def\timeofday{
\hours=\time
\minutes=\hours
\divide\hours by60
\multiply\hours by60
\advance\minutes by-\hours
\divide\hours by60
\ifnum\hours>9\else0\fi\the\hours:\ifnum\minutes>9\else
0\fi\the\minutes}
\def\predate{\date{\the\day\ \ifcase\month\or
  January\or February\or March\or April\or May\or June\or July\or
	August\or September\or October\or November\or
	   December\fi\ \the\year\ --- \timeofday\ --- Preliminary
		  Version}}

\expandafter\ifx\csname >amsstdmathmacs\endcsname\relax
\expandafter\def\csname >amsstdmathmacs\endcsname{done}
\else\endinput\fi
\def\mathcs{C^{\displaystyle *}}
\def\cs{\ifmmode\mathcs\else$\mathcs$\fi}

\def\acg{{A \rtimes_\alpha G}}
\let\bbb=\Bbb

\def\T{{\bbb T}}

\def\K{{\cal K}}
\def\iff{if and only if}
\def\Res{\operatorname{Res}}
\def\Ad{\operatorname{Ad}}
\def\Ind{\operatorname{Ind}}

\def\Prim{\operatorname{Prim}}
\def\supp{\operatornamewithlimits{supp}}
\def\ker{\operatorname{ker}}

\def\nbhd{neighborhood}
\def\id{\operatorname{id}}
\def\set#1{\{\,#1\,\}}
\let\tensor=\otimes
\def\({\bigl(}
\def\){\bigr)}
\def\restr#1{|_{{#1}}}
\def\spec#1{\specnp{(#1)}}
\def\specnp#1{{#1}^\wedge}
\newbox\hidebox
\def\spechide#1{\setbox\hidebox=\hbox{$#1$}
\hbox to\wd\hidebox{$\box\hidebox^\wedge$\hss}}
\def\dtau{{\frak t}}
\let\RR=\R

\let\ind=\Ind
\let\res=\Res
\def\ph{\varphi}

\def\Om{\Omega}
\let\bem=\em
\newsymbol\rtimes 226F
\newsymbol\subsetneqq 2324
\newsymbol\smallsetminus 2272
\let\setminus=\smallsetminus
\def\qob#1[#2]{{\cal Q}^{#1}\(#2\)}
\def\w{x}
\def\ds{{\frak s}}
\def\ttg{{\scriptstyle{\Tilde{\Tilde\Gamma}}}}
\def\ggtm{generalized Green twisting map}
\def\L{\Om}
\def\l{x}
\def\G{\Gamma}
\def\UM{\cal U}
\def\supp{\operatorname{supp}}
\def\Rep{\operatorname{Rep}}
\def\indlgh{\Ind_{\L^{N^\perp}}^{\ghat G}D}
\def\lip#1<#2,#3>{\relax_{\scriptscriptstyle{#1}}\!\!\langle #2,
#3\rangle}
\def\rip#1<#2,#3>{\langle #2, #3\rangle_{\scriptscriptstyle{#1}}}
\def\ghat{\widehat}
\def\ahat #1{\mathchoice{\widehat {#1}}{\widehat {#1}}{\hat {#1}}{\hat
{#1}}}
\def\aHat{\widehat}
\begin{document}

\title
[Crossed products]{Crossed products
whose  primitive ideal spaces are generalized trivial
$\bold{\widehat{\bold G}}$-bundles}

\author[Echterhoff]{Siegfried Echterhoff}
\address{\emergencystretch=50pt
Universit\"at-Gesamthochschule Paderborn\\
Fachbereich Mathematik-Informatik\\
Warburgerstra{\ss}e 100\\
D-4790 Paderborn\\
Germany}
\email{echter@@uni-paderborn.de}
\thanks{This research was supported by a grant from the
University of Paderborn.}

\author[Williams]{Dana P. Williams}
\address{Dartmouth College\\
Department of Mathematics\\
Bradley Hall\\
Hanover, NH 03755-3551\\
USA}
\email{dana.williams@@dartmouth.edu}
\thanks{The second author was partially supported by the National
Science Foundation.}



\date{15 August 1993}
\maketitle

\begin{abstract}
%
%
We characterize when the primitive ideal space of a crossed
product $\acg$ of a \cs-algebra $A$ by a locally compact abelian
group $G$ is a $\sigma$-trivial $\ghat G$-space for the dual $\ghat
G$-action.
Specifically, we show that $\Prim(\acg)$ is $\sigma$-trivial if
and only if the quasi-orbit space is Hausdorff, the map which
assigns to each quasi-orbit $\w$ a certain subgroup $\ttg(\alpha^\w)$
of the Connes spectrum of the system $(A_\w,G,\alpha^\w)$ is
continuous, and there is a \ggtm{} for $(A,G,\alpha)$.  Our proof
requires a substantial generalization of a theorem of Olesen and
Pedersen in which we show that there is a \ggtm{} for $(A,G,\alpha)$
\iff{} $\acg$ is isomorphic to a generalized induced algebra.
\end{abstract}

%
%
\section{Introduction}

Our Object here is to investigate the structure of the crossed
product $\acg$ of a \cs-algebra $A$ by a locally compact abelian
group $G$.
One goal of any such program is to characterize the ideal
structure via a description of the primitive ideal space
together with its Jacobson topology.  While it is theoretically
possible to describe $\Prim(\acg)$ as a set using the
Mackey-Machine as developed by Green, Gootman, and Rosenberg
\cite{green1,gr}, characterizing the topology is quite another
matter---even when $A\cong C_0(\Om,\K)$---due to the appearance of
non-trivial Mackey obstructions while employing the
Mackey-machine.  In the event that all the Mackey obstructions
vanish, considerable progress has been made in the case of
constant stabilizers \cite{pr2,rr,doir,90a,90b,90c,95a}, and
similar, but more modest, progress has been made when the
stabilizers vary continuously \cite{rw2}.
When it is not assumed that the Mackey obstructions vanish, the
situation is very murky.  Although progress has been made in
general \cite{ech1}, and especially in the case of constant
stabilizer and constant (non-vanishing) Mackey obstruction
\cite{horr,echros}, the examples in \cite{horr,echros} illustrate
the difficulties and subtleties inherent even in this restricted
approach.

One of the reasons that the approach in the constant stabilizer
(and constant Mackey obstruction) case has been so successful, is
that one can employ the machinery of classical algebraic
topology---especially the theory of principal bundles.
In the general case, the ``bundles'' that will arise will not
have isomorphic fibers, let alone be locally trivial or
principal.  However in \cite{rw2}, Iain Raeburn and the second
author introduced the notion of locally $\sigma$-trivial $G$-spaces
which we believe appropriately generalize the notion of principal
bundles.
To elaborate, recall that a locally compact $G$-space $\Om$ is called
$\sigma$-proper
if the stabilizer map $x\mapsto S_x$ is continuous and if the map
$(x,t)\mapsto (x,t\cdot x)$ defines a proper map of the quotient
$\Om\times G/\!\!\sim$ (where $(x,s)\sim(y,t)$ if and only if
$x=y$ and $ts^{-1}\in S_x$) into $\Om\times \Om$
\cite[Definition~4.1]{rw2}.  It is clear that $\sigma$-proper
spaces are the natural extensions of free and proper $G$-spaces
to non-freely acting transformation groups.  This is evidenced by
the r\^ole they and free and proper spaces play in the
investigation of crossed products with continuous trace
\cite{green2,dpw3,doir,rw2,ech3,ech2,ech4}.
Also recall that a $G$-space $\Om$ is called $\sigma$-trivial if
$\Om/G$ is Hausdorff and if it is $G$-homeomorphic to $\Om/G\times
G/\!\!\sim$ (where $(G\cdot x,s)\sim(G\cdot y,t)$ if and only
if $G\cdot x=G\cdot y$ and $ts^{-1}\in S_x$).  Of course, $\Om$
is called locally $\sigma$-trivial if $\Om$ is the union of $G$-invariant
open $\sigma$-trivial subsets \cite[Definiton~4.2]{rw2}.  The
connection between
$\sigma$-trivial spaces and
$\sigma$-properness exactly parallels that between free
and proper spaces and principal bundles: a $\sigma$-proper $G$-space
is (locally) $\sigma$-trivial if and only if there are (local)
cross sections for the orbit map \cite[Proposition~4.3]{rw2}.

A significant class of examples of $\sigma$-trivial spaces are
provided by the spectrums of transformation group \cs-algebras.
It follows from \cite{dpw4} and \cite[Theorem~5.3]{dpw2} that if
$\spec{C_0(\Om)\rtimes_\tau G}$ is Hausdorff, then it is a
$\sigma$-trivial $\ghat G$-space for the dual action.  It was
shown in \cite{rw2}, that when $A$ is non-abelian, one can have
$\spec{\acg}$ locally, but not globally, $\sigma$-trivial.
One of the original motivations for this paper was to classify
those crossed products with $\sigma$-trivial primitive ideal
space (with respect to the dual action).  Unlike \cite{rw2}, we
make {\it no} assumption on the action, the algebra $A$, or the
Mackey obstructions.

The crucial concept is that of a (generalized) Green twisting
map.  Recall that if $(A,G,\alpha)$ is an abelian dynamical
system, then a Green twisting map $\tau$ with domain $N$ is a
strictly continuous homomorphism $\tau:N\to\UM(A)$ satisfying
$\alpha_n(a)=\tau_na\tau_n^*$ and $\alpha_s(\tau_n)=\tau_n$ for
all $a\in A$, $n\in N$, and $s\in G$.
If $\tau$ is a Green twisting map for $(A,G,\alpha)$, then Olesen
and Pedersen \cite[Theorem~2.4]{op} have shown that $\acg$ is
covariantly isomorphic to an induced algebra $\Ind_{N^\perp}^G
(A\rtimes_{\alpha,\tau}G)$, where $A\rtimes_{\alpha,\tau}G$ is
Green's twisted crossed product (see \cite[\S1]{green1} for the
relevant definitions and background).
In fact, if
$\ahat\alpha$ denotes the dual $\ghat G$-action on $\acg$
and $N=\ttg(\alpha)^{\perp}$,
where $\ttg(\alpha)=\set{\gamma\in \ghat G:
\text{$\ahat\alpha_\gamma(P)=P$ for all $P\in\Prim(\acg)$}}$, then
$\ttg(\alpha)$ is a subgroup of the usual Connes spectrum
$\Gamma(\alpha)$ \cite{ped} and $A\rtimes_{\alpha,\tau} G$ is
simple if and only if $A$ is $G$-simple and $N=\{e\}$ is the
trivial subgroup \cite[Theorem~5.7]{op}.
In particular, $\Prim(\acg)$ is $\ghat G$-homeomorphic to $\ghat
N$ and is the simplest sort of $\sigma$-trivial $\ghat G$-space.
{\lineskip=0pt plus 1pt\par}
To generalize these ideas, recall that if $R:\Prim(A)\to\Om$ is a
continuous, open $G$-invariant surjection onto a locally compact
Hausdorff space $\Om$, then $\acg$ is the section algebra of a
\cs-bundle with fibers $A_x\rtimes_{\alpha^x}G$,
where $A_x=A/\ker\(R^{-1}(\set x)\)$ and $\alpha^x$ is the
induced action \cite{dpw5}.
Then, if $x\mapsto N_x$ is continuous from $\Om$
into the normal subgroups of $G$, a \ggtm{} $\dtau$ with domain
$\Om^N=\set{(x,s)\in\Om\times G:s\in N_x}$ is a ``compatible''
choice of Green twisting maps $\dtau^x$ for each quotient system
$(A_x,\, G, \,\alpha^x)$ (see Definition~\ref{deftwist} below).
Naturally, there is an associated twisted crossed product $A
\rtimes_{\alpha,\dtau}G$ which behaves as if it were fibered over
$\Om$ with fibers $A_x\rtimes_{\alpha^x,\dtau^x}G$ (see
Definition~\ref{deftcrossed} below).  Our main result
(\thmref{thm:strivial}) is that $\Prim(\acg)$ is a $\sigma$-trivial
space if and only if (1)~the quasi-orbit space $\qob G [\Prim(A)]$
is Hausdorff, (2)~the map $x\mapsto C_x=\ttg(\alpha^x)^\perp$
is continuous, and (3)~there is a \ggtm{} for $(A,\,G,\,\alpha)$ with
domain $\Om^C$.
In this event, $C_x^{\perp}=\Gamma(\alpha^x)$, the usual Connes
spectrum.
If $A$ is type~I and $(A,\,G,\,\alpha)$ is regular, then
we can replace $C$ with the symmetrizer map $\Sigma$
(\corref{cor:symm}).
Of course, as a special case of our results we
obtain conditions (\corref{cor:trivialbundle}),
which characterize when
the primitive ideal space of $A\rtimes_{\alpha}G$ is a trivial
$\ghat C$-bundle for some closed subgroup $C$ of $G$,
which generalize those given in \cite[Theorem~7.2]{90a}.

The majority of the work required to prove \thmref{thm:strivial}
involves an extensive generalization of
a theorem of Olesen and Pedersen \cite[Theorem~2.4]{op}.
Namely, we show in \thmref{theotwist} that there is a \ggtm{} $\dtau$
for $(A,\,G,\,\alpha)$ with domain $\Om^N$ if and only if the dual
system is covariantly isomorphic to a (generalized) induced system
$\(\Ind_{\Om^{N^\perp}}^{\ghat G}B,\,\ghat G,\,\ind\beta\)$ as defined
in \cite[\S3]{ech2} (or at the end of \secref{sec:2} below).
In fact, in complete analogy with \cite[Theorem~2.4]{op}, we can
take $B=A\rtimes_{\alpha,\dtau}G$.

Our work is organized as follows.
In \secref{sec:2}, we recall some of the basic definitions of
subgroup actions and subgroup crossed products from \cite{ech2}.
We also give formal definitions of \ggtm s and the resulting
twisted crossed products.  In \secref{sec:3} we prove our first
main result characterizing the crossed product of a twisted system
as an induced algebra.
Finally, in \secref{sec:4} we prove our main result
characterizing crossed products with $\sigma$-trivial primitive
ideal space.  We close with a number of interesting special
cases.

\subsection*{Acknowlegment}
The majority of the research for this article was
conducted while the second author was a visitor at the University
of Paderborn.  The second author would like to thank the first
author as well as Eberhard Kaniuth and his research group
for their warm
hospitality.
%
%
\section{Subgroup actions and generalized Green twisting maps}
\label{sec:2}

For any locally compact group $G$ we denote by $\frak K(G)$ the
space of all closed subgroups of $G$ equipped with Fell's topology
\cite{fell2}. Then $\frak K(G)$ is a compact Hausdorff space.
If $\Om$ is a locally compact space (here, and in the sequel,
locally compact always means locally compact Hausdorff)
and $H:\Om\to\frak K(G);
x\mapsto H_x$ is any continuous map, then we define
$$\Om^H=\{(x,s)\in\Om\times G: s\in H_x\},$$
which is a closed subspace of $\Om\times G$ and which may be thought as
a trivial bundle over $\Om$ with varying fibers $H_x$.

Suppose now that $A$ is a \cs-algebra such that there exists
a continuous,
open and surjective map $P:\hat{A}\to\Om$ for some locally compact space
$\Om$. Then by Lee's theorem \cite{lee2} (and
\cite[Proposition~1.6]{fell4})
we may write $A$ as the section algebra $\Gamma_0(E)$
of all continuous sections which vanish at infinity for
some \cs-bundle $P:E\to \Om$, such that the fibers $A_x$ are
isomorphic to $A/\ker P^{-1}(\{x\})$ for each $x\in \Om$.
Conversely, it was also shown by Lee that for any \cs-algebra
$A$ which is isomorphic to $\Gamma_0(E)$ for some \cs-bundle
$p:E\to\Om$ there exists a canonical continuous and open
projection from $\hat{A}$ onto $\Om$, which is just $P$ in the
situation above. We now repeat the definition of
a subgroup action as defined in \cite{ech2}.

\begin{definition}\label{defsubact}
Suppose that $A$ is a \cs-algebra such that $A$ is isomorphic to
the section algebra $\Gamma_0(E)$ for some \cs-bundle
$p:E\to\Om$. Assume further that $H:\Om\to\frak K(G);
x\mapsto H_x$ is a continuous map
such that for each $x\in \Om$ there exists a strongly continuous action
$\alpha^x$ of $H_x$ on the fiber $A_x$ such that the map
$$\Om^H\to E; (x,s)\mapsto \alpha^x_s\(a(x)\)$$
is continuous for each $a\in A$. Then $\alpha=(\alpha^x)_{x\in\Om}$
is called a {\bem subgroup action} of $\Om^H$ on $A$.
\end{definition}

If $\alpha$ is a subgroup action of $\Om^H$ on $A$, then we may
follow \cite{ech2} and form the
{\bem subgroup crossed product} $A\rtimes_{\alpha}\Om^H$ as follows:
Let $q:\Om^H\to\Om$ denote the canonical projection, and let
$\Gamma_c(q^*E)$ denote the set of continuous sections
with compact support of the pull-back bundle $q^*E$.
Note that the continuous sections of $q^*E$ are exactly given by the
continuous functions $f:\Om^H\to E$ such that $f(x,s)\in A_x$ for all
$x\in \Om$. We define convolution, involution and norm on
$\Gamma_c(q^*E)$ by the rules
\begin{gather*}
f*g(x,s)=\int_{H_x}f(x,t)\alpha^x_t\(g(x,t^{-1}s)\)\,d_{H_x}t,\\
f^*(x,s)=\Delta_{H_x}(s^{-1})\alpha^x_s(f(x,s^{-1})^*),\\
\intertext{and}
\|f\|_1=\sup_{x\in\Om}\int_{H_x}\|f(x,s)\|\,d_{H_x}s,
\end{gather*}
where $\Delta_{H_x}$ denotes the modular function on $H_x$ and
$(d_H)_{H\in\frak K(G)}$ is a smooth choice of Haar measures on
$\frak K(G)$ \cite[p.~908]{gl2}.
If $L^1(\Om^H,A,\alpha)$ denotes the completion of $\Gamma_c(q^*E)$ with
respect to $\|\cdot\|_1$, then the subgroup crossed product
$A\rtimes_{\alpha}\Om^H$ is defined as the enveloping \cs-algebra of
$L^1(\Om^H,A,\alpha)$.

It has been shown in \cite{ech2} that  $\spec{A\rtimes_{\alpha}\Om^H}$
can be identified with the set of pairs
$$\set{(x, \rho\times V): x\in\Om, \rho\times V\in
\spec{A_x\rtimes_{\alpha^x}H_x}}$$
by defining
$$(x,\rho\times V)(f)=\rho\times V\(f(x,\cdot)\).$$
More generally. Every representation of a fiber $A_x\rtimes_{\alpha^x}H_x$
defines a representation of $A\rtimes_{\alpha}\Om^H$ in this way
and the collection $\cal S(A\rtimes_{\alpha}\Om^H)$ of
all equivalence classes of
those representation is called the space of {\bem subgroup representations
of $A\rtimes_{\alpha}\Om^H$}. If we restrict to representations with
dimension bounded by a fixed cardinal, say $\aleph$, then we may topologize
$\cal S(A\rtimes_{\alpha}\Om^H)$ by viewing $\cal S(A\rtimes_{\alpha}\Om^H)$
as a subset of the space $\Rep(A\rtimes_{\alpha}\Om^H)$ of all equivalence
classes of representations of $A\rtimes_{\alpha}\Om^H$
with dimension bounded by $\aleph$
equipped with Fell's
inner hull kernel topology (see \cite{fell6,fell3}).
Note that the canonical projection of $\cal S(A\rtimes_{\alpha}\Om^H)$
onto $\Om$ is always continuous. Moreover, it has been shown in
\cite[Corollary~4]{ech2}
that in case where $H_x$ is amenable for all $x\in \Om$,
the restriction of this projection to $\spec{A\rtimes_{\alpha}\Om^H}$
is also open, so that by Lee's theorem we may
write $A\rtimes_{\alpha}\Om^H$ as a section algebra $\Gamma_0(F)$
of a section bundle $r:F\to \Om$ with fibers $A_x\rtimes_{\alpha^x}H_x$.

The definition of a subgroup action on a \cs-algebra as above was
motivated by the desire to be able to restrict an action of $G$ on
a \cs-algebra $A$ to a set of continuously varying subgroups of $G$.
Before we make this more precise we first want to introduce the
notion of a regularization of a dynamical system.

\begin{definition}\label{defreg}
Let $(A,\,G,\,\alpha)$ be a dynamical system. A {\bem regularization}
$(\Om,R)$ of $(A,\,G,\,\alpha)$ consists of a locally
compact
$G$-space $\Om$ together with
a $G$-equivariant continuous map $R:\Prim(A)\to\Om$. If, in addition, $R$  is
open and surjective, then $(\Om,R)$ is called an {\bem open regularization}
of $(A,\,G,\,\alpha)$.  If the action of $G$ on $\Om$ is trivial,
then $(\Om,R)$ is called a {\bem $G$-invariant regularization}.
\end{definition}
\begin{remark}
Since $\Om$ is Hausdorff, any continuous map from $\hat A$ into
$\Om$ must factor through $\Prim(A)$.  Thus whenever convenient
we can, and will, think of the regularizing maps $R$ as being
defined on $\hat A$.
\end{remark}

Note that the definition of a regularization as given above is much weaker
than the definition of regularizations as defined in \cite{ech2}.
Please see \cite[\S1]{ech2} for situations where both
definitions automatically coincide. We now come to the motivating example
for subgroup actions.

\begin{ex}\label{ex1}
Suppose  that $(\Om,R)$ is an open regularization of
$(A,\,G,\,\alpha)$. Assume that
$H:\Om\to \frak K(G)$ is a continuous map such that $H_x\subseteq S_x$
for all $x\in\Om$, where $S_x$ denotes the stabilizer of $x$.
Then we may define for each $x\in\Om$ an action $\alpha^x$ of $H_x$
on the fiber $A_x=A/\ker R^{-1}(\{x\})$ of the corresponding \cs-bundle
$r:E\to\Om$ with $A\cong \Gamma_0(E)$,
by first restricting $\alpha$ to $H_x$, and then taking the
canonical action of $H_x$ on the quotient of $A$ by the $H_x$-invariant
ideal $\ker R^{-1}(\{x\})$. Then $(\alpha^x)_{x\in\Om}$ is a subgroup
action of $\Om^H$ on $A$. We will call this action the
{\bem restriction of $\alpha$ to $\Om^H$} and will denote it usually
also by the letter $\alpha$.
This will cause no confusion as it will be clear from context
whether
$\alpha$ refers to an action
of $G$ on $A$, or the restriction of an action to $\Om^H$ as
defined above.
\end{ex}

The construction of the subgroup action in the example above
gives rise to many different constructions of subgroup algebras.
One very important case is the following.

\begin{ex}\label{ex2}
Suppose that $(A,\,G,\,\alpha)$ is a dynamical system,
$\Om$ is a locally compact space and
$H:\Om\to \frak K(G)$ is a continuous map.
Let $\id\otimes\alpha$ denote the diagonal action of $G$ on $C_0(\Om,A)$
with respect to the trivial action of $G$ on $\Om$ and let
$R:\specnp{C_0(\Om,A)}\to\Om$ denote the canonical projection.
Then $(\Om,R)$ becomes an open regularization for
$(C_0(\Om,A),\,G,\,\id\otimes\alpha)$ and
as in Example \ref{ex1} we may form the subgroup action
$\id\otimes\alpha$ of $\Om^H$ on $C_0(\Om,A)$
and the subgroup crossed product
$C_0(\Om,A)\rtimes_{\id\otimes\alpha}\Om^H$, which we will denote in
the following simply by $\cs (\Om^H,A,\alpha)$.

Note that we could have constructed $\cs (\Om^H,A,\alpha)$ easily by defining
convolution and involution on $C_c(\Om^H,A)$
by the usual operations on the fibers $C_c(H_x,A)$ with respect to the
restriction of $\alpha$ to $H_x$.
The space $\cal S(\Om^H,A,\alpha)$
of subgroup representations of $\cs (\Om^H,A,\alpha)$ consists of
all pairs $(x, \rho\times V)$ with $x\in \Om$ and
$\rho\times V\in\Rep(A\rtimes_{\alpha}H_x)$.
In case where $\Om$ is equal to $\frak K(G)$ and $H$ equals the identity
$I:\frak K(G)\to \frak K(G)$, we obtain the subgroup algebra
$\cs (\frak K(G)^I,A,\alpha)$ whose subgroup representations give the
collection of {\bem all} subgroup representations for the system
$(A,\,G,\,\alpha)$.
\end{ex}

The following proposition shows that there is a very strong connection between
the constructions in Examples \ref{ex1} and \ref{ex2}.

\begin{prop}\label{propisom1}
Suppose that $(\Om,R)$ is an open regularization of $(A,\,G,\,\alpha)$ and
assume that $H:\Om\to\frak K(G)$
is a continuous map such that $H_x$
is contained in the stabilizer $S_x$ for all $x\in \Om$.
Let $A\rtimes_{\alpha}\Om^H$ be the subgroup crossed product of $A$
by $\Om^H$ with respect to the restriction of $\alpha$ to $\Om^H$,
and let $\cs (\Om^H,A,\alpha)$ denote the subgroup algebra as constructed
in Example \ref{ex2}.
Let $\cal S_R^H$ denote the set of all irreducible representations
of $\cs (\Om^H,A,\alpha)$ which are given by pairs $(x,\rho\times V)$
such that $\rho\times V\in \spec{A\rtimes_{\alpha}H_x}$ satisfies
$\ker\rho\supseteq\ker R^{-1}(\{x\})$ (i.e.\ $\rho\in\Rep(A_x)$).
Then $\cs (\Om^H,A,\alpha)/\ker\cal S_R^H$
is canonically isomorphic to $A\rtimes_{\alpha}\Om^H$.
\end{prop}
\begin{pf}
Let $r:E\to\Om$ be the \cs-bundle given by $R:\hat{A}\to\Om$ and Lee's
theorem, and let as usual $q^*E$ denote the pull-back of $E$ via the
canonical projection $q:\Om^H\to\Om$. Then define
$$\Phi:C_c(\Om^H,A)\to\Gamma_c(q^*E); \Phi(f)(x,s)= f(x,s)(x).$$
Then it is easily seen  that $\Phi$ extends to
a surjective $*$-homomorphism from $\cs (\Om^H,A,\alpha)$ onto
$A\rtimes_{\alpha}\Om^H$ with kernel $\ker \cal S_R^H$ (compare with
\cite[Proposition~9]{ech2}).
\end{pf}

Recall that a Green twisting
map $\tau$ for a dynamical system $(A,\,G,\,\alpha)$
is a strictly continuous homomorphism $\tau:N_{\tau}\to\cal U(A)$
from a closed normal subgroup $N_{\tau}$ of $G$ into the group of unitaries
$\cal U(A)$ in the multiplier algebra $\cal M(A)$ of $A$ such that
the conditions
$$\alpha_n(a)=\tau_na\tau_n^*\;\;\;\text{and}\;\;\;
\alpha_s(\tau_n)=\tau_{sns^{-1}}$$
are satisfied for all $n\in N_{\tau}$, $a\in A$ and $s\in G$.
The twisted crossed product $A\rtimes_{\alpha,\tau}G$
is defined as the quotient of $A\rtimes_{\alpha}G$ by the intersection of
all kernels of representations $\pi\times U$ of $A\rtimes_{\alpha}G$
which preserve $\tau$ in the sense that
$$\pi(\tau_n)=U_n\;\;\;\text{for all}\;\;\;n\in N_{\tau}.$$

Our first aim here is to generalize the notion of Green twisting
maps to
something which has as the domain a set of varying subgroups of $G$,
rather than just  one fixed normal subgroup $N_{\tau}$ of $G$.
For this note first that if $(\Om,R)$ is a $G$-invariant
open regularization of a system $(A,\,G,\,\alpha)$, then
$\ker R^{-1}(\{x\})$ is
a $G$-invariant ideal in $A$ for each $x\in \Om$ and
we have
a canonical action $\alpha^x$ of $G$ on each fiber
$A_x=A/\ker R^{-1}(\{x\})$.

\begin{definition}\label{deftwist}
Assume that $(\Om,R)$ is a $G$-invariant open regularization of the
dynamical system $(A,\,G,\,\alpha)$, $r:E\to\Om$ the corresponding
\cs-bundle with fibers $A_x=A/\ker R^{-1}(\{x\})$,
and let $N:\Om\to\frak K(G)$ be a
continuous map such that $N_x$ is normal in $G$ for all $x\in \Om$.
A {\bem \ggtm} with domain $\Om^N$ for
$(A,\,G,\,\alpha)$ is a map
$\dtau:\Om^N\to \bigcup_{x\in\Om}\cal U(A_x)$ such that
\begin{enumerate}
\item Each map $\dtau^x:N_x\to \cal U(A_x); \dtau^x_n=\dtau(x,n)$
is a Green twisting map for the system $(A_x,\,G,\,\alpha^x)$.
\item The maps $\Om^N\to E$ given by
$(x,n)\mapsto \dtau_n^x a(x)$ and $(x,n)\mapsto
a(x)\dtau^x_n$ are continuous for each
$a\in A=\Gamma_0(E)$.
\end{enumerate}
If $\dtau$ is a \ggtm{} for $(A,\,G,\,\alpha)$,
then we call $(A,\,G,\,\alpha,\,\dtau)$ a {\bem generalized twisted
dynamical system}.
\end{definition}
\begin{remark}\label{rem:2.8}
Note that if we do not require that $\alpha^x_s(\dtau^x_n) =
\dtau^x_n$ for each $x\in\Om$, $s\in G$ and $n\in N_x$---that is,
each $\dtau^x$ is merely a strictly continuous
homomorphism into $\cal U(A_x)$ implementing $\alpha^x$---then we
simply have a {\bem unitary\/} subgroup action of $\Om^N$ on $A$
as defined in \cite[Definition~5.1]{rw2}.
\end{remark}

As in the case of ordinary twisted dynamical systems, there is a
canonical procedure for defining twisted crossed products of
generalized twisted dynamical systems. The easiest definition is the
following.

\begin{definition}\label{deftcrossed}
Suppose that $\dtau$ is a \ggtm{}
for $(A,\,G,\,\alpha)$
with domain $\Om^N$. We say that a representation $\pi\times U$ of
$(A,\,G,\,\alpha)$ {\bem preserves $\dtau$} if $\pi(\dtau^x_s)=U_s$ for all
$(x,s)\in\Om^N$ such that $\ker\pi\supseteq\ker R^{-1}(\{x\})$.
Let $I_{\dtau}$ be the intersection of all kernels of representations
of $A\rtimes_{\alpha}G$ which preserve $\dtau$.
Then we define the {\bem generalized twisted crossed product}
$A\rtimes_{\alpha,\dtau}G$ as the quotient of $A\rtimes_{\alpha}G$
by $I_{\dtau}$.
\end{definition}

The following observation will be useful in considering the structure of
generalized twisted crossed products.  The proof is
straightforward.
\begin{lem}\label{lemid}
Suppose that $(\Om,R)$ is a $G$-invariant open regularization of
$(A,\,G,\,\alpha)$. We identify $G$ with the constant map
$G:\Om\to\frak K(G); x\mapsto G$. Let $\alpha=(\alpha^x)_{x\in\Om}$
be the subgroup action of $\Om^G$ on $A$ as in Example \ref{ex1}.
Then
$$\Psi:C_c(G,A)\to \Gamma_c(q^*E); \Psi(f)(x,s)=f(s)(x)$$
extends to an isomorphism between $A\rtimes_{\alpha}G$ and
$A\rtimes_{\alpha}\Om^G$.
\end{lem}

Note that by the identification of $A\rtimes_{\alpha}G$ with
$A\rtimes_{\alpha}\Om^G$ via \lemref{lemid}, we
know that for every irreducible representation $\pi\times U$
of $A\rtimes_{\alpha}G$ there exists exactly one $x\in\Om$
such that $\ker\pi\supseteq \ker R^{-1}(\{x\})$.
Hence it follows that the collection of irreducible representations
of $A\rtimes_{\alpha}G$ which preserve $\dtau$ may be viewed as
the union over all $x\in\Om$ of  the spaces
$\spec{A_x\rtimes_{\alpha^x,\dtau^x}G}$,
where $A_x\rtimes_{\alpha^x,\dtau^x}G$ denotes the usual
twisted crossed product with respect to the ordinary twisted action
$(\alpha^x,\dtau^x)$ of $G$ on $A_x$.
It is straight-forward to see that $I_{\dtau}$ is equal to the intersection
of all kernels of the irreducible representations which preserve $\dtau$.
We will later see that
these are exactly the irreducible representations of
$A\rtimes_{\alpha,\dtau}G$.

As for ordinary twisted crossed products it is quite useful to have another
realization of the generalized twisted crossed product. For this
we have to define first what we will understand as a quotient of
$G$ by $\Om^N$.

\begin{definition}\label{defquot}
Assume that $G$ is a locally compact group, $\Om$ a locally compact space
and $H:\Om\to \frak K(G)$ is a continuous map. Then we define an
equivalence relation $\sim_H$ on $\Om\times G$ by
$$(x,s)\sim_H(y,t)\text{ if and only if
$x=y$ and $s\in tH_x$}$$
The quotient space $\Om\times_HG=(\Om\times G)/\!\!\sim_H$
under this equivalence
relation is called the {\bem quotient of $G$ by $\Om^H$}.
\end{definition}

Note that the elements of $\Om\times_HG$ as in the definition above are just
given by the pairs $(x,\dot{s})$, where $\dot{s}$ denotes the left
coset $sH_x$ of $s\in G$.

Assume now that $N:\Om\to \frak K(G)$ is a continuous map such that $N_x$
is normal in $G$ for all $x\in\Om$. Then we may define
(left) Haar measures
on $G/N_x$ by normalizing with respect to a fixed Haar measure on $G$
and a smooth choice of
Haar measures on $\frak K(G)$ such that
$$\int_G g(s)\, ds=\int_{G/N_x} \int_{N_x}g(sn)\,d_{N_x}n\,d_{G/N_x}\dot{s},$$
for any $g\in C_c(G)$.
(While a similar family of measures in defined in \cite[Lemma~2.20]{dpw3},
the reader should be cautioned that in the case the subgroups $S_x$
are normal the measures $\mu_x$ defined in \cite{dpw3} are {\it
right} Haar measures (see \cite[Lemma~2.21]{dpw3}).)
We record some properties that these measures enjoy for later
use.

\begin{lem}\label{lemint}
Suppose we have chosen Haar measures on the quotient groups $G/N_x$
as above.
\begin{enumerate}
\item
Let $\Lambda$ be a locally compact space and suppose $F\in
C_c\(\Lambda,(\Om\times_NG)\)$.  Then
\begin{equation*}
f(\lambda,x)=\int_{G/N_x} F(\lambda,x,\dot s)\,d_{G/N_x}\dot s
\end{equation*}
defines an element of $C_c(\Lambda\times\Om)$.
\item
For each $f\in C_c(\Om\times_NG)$,
\begin{equation*}
\tilde f(x)=\int_{G/N_x} f(x,\dot s)\,d_{G/N_x}\dot s
\end{equation*}
defines and element of $C_c(\Om)$.
\item
The map $(x,\dot s)\mapsto \Delta_{G/N_x}(\dot s)$ is continuous
on $\Om\times_NG$.
\end{enumerate}
\end{lem}
\begin{pf}
The first assertion can be proved along the same lines as
\cite[Lemmas 2.5(iv) and 2.22]{dpw3}.
The second and third assertions follow from the first: (1)
follows immediately and (2) follows by considering $F(\dot r,x,
\dot s)=f(x,\dot r\dot s)$ for appropriate choices of $f$ in $C_c
(\Om\times_NG)$.
\end{pf}

If $\dtau$ is a generalized twist for
$(A,\,G,\,\alpha)$ with domain $\Om^N$ as in Definition \ref{deftwist},
then we define $C_c(G,A,\dtau)$ as the set of all continuous $A$-valued
functions $f$ on $G$ which satisfy
$$f(ns)(x)=f(s)(x)\dtau^x_{n^{-1}}$$
for all $s\in G$ and $(x,n)\in \Om^N$ such that the resulting function
$$(x,\dot{s})\mapsto\|f(s)(x)\|$$
has compact support in $\Om\times_NG$.
Then, as usual, we define multiplication, involution and norm on
$C_c(G,A,\dtau)$ by
\begin{gather}
f*g(s)(x)=\int_{G/N_x}f(t)(x)\alpha^x_t\(g(t^{-1}s)(x)\)\,d_{G/N_x}t,
\label{eq:star1}\\
\label{eq:star2}
f^*(s)(x)=\Delta_{G/N_x}(s^{-1})\alpha^x_s(f(s^{-1})(x)^*), \\
\intertext{and}
\notag
\|f\|_1=\sup_{x\in \Om}\int_{G/N_x} \|f(t)(x)\|\,d_{G/N_x}t.
\end{gather}

Note that if we define $f^x$ by $f^x(s)= f(s)(x)$ for all $x\in\Om$
and $f\in C_c(G,A,\dtau)$, then $f^x$ is an element of the
dense subspace $C_c(G,A_x,\dtau^x)$ of the ordinary twisted crossed product
$A_x\rtimes_{\alpha^x,\dtau^x}G$ as constructed by Green in
\cite[Section 1]{green1}. Moreover, if we define $\|\cdot\|_1^x$
on $C_c(G,A_x,\dtau^x)$ by
$$\|g\|_1^x=\int_{G/N_x} \|g(s)\|\,d_{G/N_x}\dot{s},$$
then it was shown also in \cite[Section 1]{green1} that
the twisted crossed product $A_x\rtimes_{\alpha^x,\dtau^x}G$ is just
the enveloping \cs-algebra of the completion $L^1(G,A_x,\dtau^x)$
of $C_c(G,A,\dtau)$ with respect to $\|\cdot\|_1^x$.
Now it follows directly from Lemma \ref{lemint} that for each
$f\in C_c(G,A,\dtau)$ the map $\Om\to \RR; x\mapsto\|f^x\|^x_1$ is continuous
with compact support. Thus we see that the completion $L^1(G,A,\dtau)$
of $C_c(G,A,\tau)$ with respect to $\|\cdot\|_1$ is
isomorphic to $\Gamma_0(F)$, for some Banach bundle $p:F\to\Om$
with fibers $B_x=L^1(G,A_x,\dtau^x)$. Observe also that multiplication
and involution on $C_c(G,A,\dtau)$ are just defined fiberwise by
the multiplication and involution on the fibers $L^1(G,A_x,\dtau^x)$.
We use this observation to prove

\begin{prop}\label{proptwistL1}
Multiplication and involution on $C_c(G,A,\dtau)$, as defined in
Equations \eqref{eq:star1}~and \eqref{eq:star2},
are well defined and extend to all of
$L^1(G,A,\dtau)$. Thus $L^1(G,A,\dtau)$ is the section algebra
$\Gamma_0(F)$ of a Banach $*$-algebra
bundle $p:F\to\Om$ with fibers $L^1(G,A_x,\dtau^x)$.
Moreover, the map
$$\Phi:C_c(G,A)\to C_c(G,A,\tau);
\Phi(f)(s)(x)=\int_{N_x}f(sn)(x)\dtau^x_{sns^{-1}}\,d_{N_x}n$$
extends to a $*$-homomorphism from $A\rtimes_{\alpha}G$ onto
the enveloping \cs-algebra $\cs (G,A,\dtau)$ of $L^1(G,A,\dtau)$ with
kernel $I_{\dtau}$. Thus $\cs (G,A,\dtau)$ is canonically isomorphic to
$A\rtimes_{\alpha,\dtau}G$. Moreover, the irreducible representations
of $A\rtimes_{\alpha,\dtau}G$ are exactly the irreducible representations
of $A\rtimes_{\alpha}G$ which preserve $\dtau$.
\end{prop}
\begin{pf}
That \eqref{eq:star1}~and \eqref{eq:star2}
define elements of $C_c(G,A,\dtau)$ is clear except for possibly
the continuity of $f*g$.  But if $f,g\in C_c(G,A,\dtau)$, then
the function $(x,r)\mapsto f(r)(x)\alpha^x_r\(g(r^{-1}s)(x)\)$
defines an element $\phi_s\in C_c(\Om\times_N G,A)$.  In fact if
$\rho:\Om\times G\to \Om\times_NG$ is the quotient map, and if
$s_0\in G$ is fixed, there are compact sets $C\subseteq	\Om$
and $K\subseteq G$ such that $\supp(\phi_s) \subseteq
\rho(C\times K)$ for all $s$ near $s_0$ \cite[Lemma~2.3]{dpw3}.
Now choose a generalized Bruhat approximate cross section
(\cite[Proposition~2.18]{dpw3}) $b\in C_c(\Om\times G)$ such that
\begin{equation*}
\int_{N_x}b(x,st)\,d_{N_x}t = 1
\end{equation*}
provided $x\in C$ and $s\in KN_x$.  Then  we can define a
continuous function $F:G\to A$ by
\begin{equation} \label{eq:func}
F(s)(x)=\int_G \phi_s(x,\dot r)b(x,r)\,dr.
\end{equation}
Of course, near $s_0$, \eqref{eq:func} equals
\begin{equation*}
\int_{G/N_x}f(r)(x)\alpha_r^x\(g(r^{-1}s)(x)\)d_{G/N_x}\dot r =
f*g(s)(x).
\end{equation*}

Next
observe that by the properties of a continuous choice of
Haar measures, together with the continuity of $\dtau$, it is straight-forward
to see that $\Phi(f)\in C_c(G,A,\dtau)$ for any $f\in C_c(G,A)$.
In fact, $\Phi$ maps $C_c(G,A)$ onto $C_c(G,A,\dtau)$.  To see
this, consider $f\in C_c(G,A,\dtau)$ with support contained in
$\rho(C\times K)$ as above.  Also as above choose $b\in C_c(\Om
\times G)$ with
\begin{equation*}
\int _{N_x}b(x,st)\,d_{N_x}t = 1
\end{equation*}
for all $x\in C$ and $s\in KN_x$.  Then $F\in C_c(G,A)$ defined
by $F(s)(x)=f(s)(x)b(x,s)$ satisfies $\Phi(F)=f$.

Now, for each $x\in \Om$, let us look at the map
$$\Phi^x:C_c(G,A)\to L^1(G,A_x,\dtau^x); \Phi^x(f)=\Phi(f)^x=
\int_{N_x}f(sn)(x)\dtau^x_{sns^{-1}}\,d_{N_x}n$$
Observe that $\Phi^x$ is the composition of the natural maps
$\Psi^x:C_c(G,A)\to C_c(G,A_x)$ given by $\Psi^x(f)(s)=f(s)(x)$
and $\ph^x:C_c(G,A_x)\to L^1(G,A_x,\dtau^x)$ which is given by
integration against $\dtau^x$. Since all these maps extend to $*$-homomorphisms
from the appropriate $L^1$-completions onto their images, we conclude that
each $\Phi^x$ extends to an $*$-homomorphism from
$L^1(G,A)$ onto $L^1(G,A_x,\dtau^x)$.
Thus we see that $\Phi$ is a $*$-homomorphism with respect to
multiplication and involution on $C_c(G,A,\dtau)$ as defined above.

It has now become clear that $L^1(G,A,\dtau)$ is equal to $\Gamma_0(F)$
for a Banach $*$-algebra bundle $p:F\to\Om$ with fibers $L^1(G,A_x,\dtau^x)$.
But for these algebras it is clear that the irreducible
representations of $L^1(G,A,\dtau)$ are canonically
given by the representations of the fibers $L^1(G,A_x,\dtau^x)$.
This shows that the representations of $L^1(G,A,\dtau)$, and hence
also of $\cs (G,A,\dtau)$ are given by the union over all
$x\in\Om$ of the spaces $\spec{A_x\rtimes_{\alpha^x,\dtau^x}G}$.

Finally, if $f\in C_c(G,A)$ and if $\pi\times U\in
\spec{A_x\rtimes_{\alpha^x,\dtau^x}G}$ is viewed
both as an element
of $\spec{A\rtimes_{\alpha}G}$ on the one hand,
and as an element of
$\specnp{\cs (G,A,\dtau)}$ on the other,
then it follows directly from the constructions
that $\pi\times U\(\Phi(f)\)=\pi\times U(f)$ for all $f\in C_c(G,A)$.
Thus $\Phi$ extends to a $*$-homomorphism from $A\rtimes_{\alpha}G$
onto $\cs (G,A,\dtau)$ such that $\ker\Phi$ is equal to $I_{\dtau}$.
But this implies that $\Phi$ factors through an isomorphism between
$A\rtimes_{\alpha,\dtau}G$ and $\cs (G,A,\dtau)$, which finishes the proof.
\end{pf}

As in the case of ordinary Green twisting maps there is a canonical
example for generalized Green twisting maps.

\begin{ex}\label{ex3}
Suppose that $(\Om,R)$ is a $G$-invariant open regularization of
$(A,\,G,\,\alpha)$, and let $N:\Om\to \frak K(G)$ be a continuous map
such that $N_x$ is a normal subgroup of $G$ for all $x\in\Om$.
Then we may restrict $\alpha$ to $\Om^N$ as in Example \ref{ex1}
and we may form the subgroup crossed product $A\rtimes_{\alpha}\Om^N$.
Let $\gamma^N$ denote the action of $G$ on $A\rtimes_{\alpha}\Om^N$,
which is defined by
\begin{equation}
\label{eq:gammaN}
\gamma^N_s(f)(x,n)= \delta(x,s)\alpha^x_s\(f(x,n^{-1}s)\),
\end{equation}
$f\in \Gamma_0(q^*E)$, where $\delta:\Om\times G\to\RR^+$ is defined by
$\delta(x,s)=\Delta_{G/N_x}(\dot s)\Delta_G(s^{-1})$.
(Recall that $\delta$ is continuous by virtue of \lemref{lemint}.)

Note that the quotient action $\gamma^{N_x}$ of $G$ on
$A_x\rtimes_{\alpha^x}N_x$ is the action defined by Green in
\cite[Proposition~1]{green1}.  Consequently, there is a Green
twisting map $\dtau^{N_x}:N_x\to \cal
U(A_x\rtimes_{\alpha^x}N_x)$ for the quotient system
$(A_x\rtimes_{\alpha^x}N_x, \,G,\,\gamma^{N_x})$.

Recall that $A\rtimes_{\alpha}\Om^N$ is
a section algebra $\Gamma_0(D)$ of a \cs-bundle $p:D\to\Om$ with
fibers $A_x\rtimes_{\alpha^x}N_x$ \iff{}
the canonical projection $P:\spec{A\rtimes_{\alpha}\Om^N}\to\Om$ is open
(recall that this is always true if all $N_x$ are amenable).
If $P$ is open,
then we claim that the map
$$\dtau^N:\Om^N\to \bigcup_{x\in\Om}\cal U(A_x\rtimes_{\alpha^x}N_x);
\dtau^N(x,n)=\dtau^{N_x}_n$$
is a generalized Green twisting map for the system
$(A\rtimes_{\alpha}\Om^N,\,G,\,\gamma^N)$.
In order to see this we only have to check the continuity of the maps
from
$\Om^N\to D$ defined, for each $d\in A\rtimes_{\alpha}\Om^N$, by
\begin{equation*}
(x,n)\mapsto \dtau^{N_x}_nd(x)\quad\text{and}\quad
(x,n)\mapsto d(x)\dtau^{N_x}_n.
\end{equation*}
But if $d=f\cdot a$
is defined as $f\cdot a(x,s)=f(x,s)a(x)$ for
$f\in C_c(\Om^N)$ and $ a\in A$,
we see that $\(\dtau^x_n f\cdot a(x)\)(s)=f(x,n^{-1}s)\alpha^x_n\(a(x)\)$
and $(f\cdot a(x)\dtau^x_n)(s)=\Delta_{N_x}(n^{-1})f(xsn^{-1})a(x)$,
and it is clear that these expressions are continuous with respect to
all variables $x$, $s$ and $n$ (see \cite[Lemma~2.5(ii)]{dpw3}).
Thus the two maps given above are easily seen
to be continuous for such $d$. Since the set
$\set{f\cdot a:f\in C_c(\Om^N), a\in A}$ generates a dense subset of
$A\rtimes_{\alpha}\Om^N$, the continuity follows for all $d\in
A\rtimes_{\alpha}\Om^N$.

Next, let us consider the generalized twisted crossed product
$(A\rtimes_{\alpha}\Om^N)\rtimes_{\gamma^N,\dtau^N}G$. In analogy to
the situation of ordinary twisted crossed products
\cite[Proposition~1]{green1}, we want to show that
$(A\rtimes_{\alpha}\Om^N)\rtimes_{\gamma^N,\dtau^N}G$ is isomorphic to
$A\rtimes_{\alpha}G$.
For this let $r:E\to \Om$ denote the \cs-bundle belonging to $A$ and let
$\Gamma_c\(q^*(E)\)$ denote the dense subalgebra of $A\rtimes_{\alpha}\Om^N$
of continuous sections with compact support of the pull-back bundle
$q^*E$ with respect to the projection $q:\Om^N\to\Om$.
We define
$\Psi:C_c(G,A)\to C_c(G,\Gamma_c(q^*E), \dtau^N)$ by
$$\(\Psi(f)(s)\)(x,n)=\delta(x,s)f(ns)(x).$$
Then it follows
from the remarks preceding \cite[Proposition 8]{ech1} that
$\Psi$ is given fiberwise by the canonical isomorphism between the fibers
$A_x\rtimes_{\alpha^x}G$ and
$(A_x\rtimes_{\alpha^x}N_x)\rtimes_{\gamma^{N_x},\dtau^{N_x}}G$.
Since it is also clear that $\Psi\(C_c(G,A)\)$ is invariant under
pointwise multiplication with bounded continuous functions on $\Om$
we conclude that
$\Psi\(C_c(G,A)\)$ is dense in
$(A\rtimes_{\alpha}\Om^N)\rtimes_{\gamma^N,\dtau^N}G$.
Thus $\Psi$ extends to the desired isomorphism
between $A\rtimes_{\alpha}G$ and
$(A\rtimes_{\alpha}\Om^N)\rtimes_{\gamma^N,\dtau^N}G$.
\end{ex}

Before we close this section we want to recall from \cite{ech2}
some basic results about $\sigma$-trivial regularizations and
induced dynamical systems.
For this recall first that a locally compact $G$-space $\Om$ is called
a {\bem $\sigma$-proper $G$-space} if the stabilizer map
$S:\Om\to\frak K(G);x\mapsto S_x$ is continuous and the canonical map
$$\Om\times_SG\to\Om\times\Om; (x,\dot{s})\mapsto(x,sx)$$
is proper in the usual sense that inverse images of compact sets are compact.
If $\Om$ is a $\sigma$-proper $G$-space, then we say that $\Om$ is
{\bem $\sigma$-trivial}, if there exists a continuous section
$\frak s:\Om/G\to\Om$ for the quotient map $\Om\to\Om/G$.
Moreover, $\Om$ is called {\bem locally $\sigma$-trivial}, if $\Om$
is a $\sigma$-proper $G$-space with local sections, which just means
that each $x\in \Om$ has a $G$-invariant neighborhood $U$ such that
$U$ is a $\sigma$-trivial $G$-space.

Now let $(A,\,G,\,\alpha)$ be any dynamical system, and let $(\Om,R)$
be a regularization of $(A,\,G,\,\alpha)$. Then we say that
$(\Om,R)$  is a {\bem $\sigma$-trivial regularization} of $(A,\,G,\,\alpha)$
if $\Om$ is a $\sigma$-trivial $G$-space.

It was shown in \cite{ech2}
that the dynamical systems which possess a $\sigma$-trivial {\bem open}
regularization are exactly those which are induced from
a subgroup action as follows.
Start with any action $\beta$ of $\Om^H$ on a \cs-algebra
$B$ as in Definition \ref{defsubact}, where $B=\Gamma_0(E)$ for some
\cs-bundle $p:E\to\Om$, and $H:\Om\to \frak K(G)$ is a continuous map.
Then the {\bem induced \cs-algebra} $\Ind_{\Om^H}^GD$ is defined
in \cite[Definition~5]{ech2}
as
\begin{multline*}
\Ind_{\Om^H}^GB=\{\,f\in C_b(G,B): \text{$f(sh^{-1})(x)
=\alpha^x_h\(f(s)(x)\)$
for all $(x,h)\in \Om^H$,} \\
\text{and such that $
(x,\dot{s})\to \|f(s)(x)\|$ is a function in
$ C_c(\Om\times_HG)$}\,\}.
\end{multline*}
The induced action $\Ind\beta$ of $G$ on $\Ind_{\Om^H}^GB$ is given by
$\(\Ind\beta_s(f)\)(t)=f(s^{-1}t)$
for all $s,t\in G$.

For each $\pi\in\hat{B}$ and $s\in G$ let $M(\pi,s)$ denote the
irreducible representation of $\Ind_{\Om^H}^GB$ given by
$$M(\pi,s)(f)=\pi\(f(s)\),\qquad f\in\Ind_{\Om^H}^GB.$$
Then $M(\pi,s)=M(\rho,t)$
in $\spec{\Ind_{\Om^H}^GB}$ if and only if
$(\pi,s)\sim(\rho,t)$, where $\sim$ denotes the
equivalence relation on $\hat{B}\times G$ which is given by
\begin{multline*}
\text{$(\pi, s)\sim(\rho, t)$ if and only if $P(\pi)=P(\rho)=x$
and there exists} \\
\text{$h\in H_x$
such that $(\pi,s)=(\rho\circ\beta^x_h, th)$.}
\end{multline*}
Here, $P:\hat{B}\to\Om$ denotes the canonical projection corresponding
to $p:E\to\Om$.
Moreover, each representation of $\Ind_{\Om^H}^GB$ is given in this way
and $\spec{\Ind_{\Om^H}^GB}$ is homeomorphic
to $(\hat{B}\times G)/\!\!\sim$
(see \cite[Proposition 10]{ech2}).

Hence we observe easily that there
exists a canonical continuous open and surjective $G$-map
$R:\spec{\Ind_{\Om^H}^GB}\to \Om\times_HG$, which is given by mapping
$M(\pi,s)\in
\spec{\Ind_{\Om^H}^GB}$ to the pair $(P(\pi),\dot{s})\in
\Om\times_HG$.  Thus the pair
$(\Om\times_HG,R)$ becomes an open $\sigma$-trivial regularization
of $(\Ind_{\Om^H}^GB,\,G,\,\Ind\beta)$.

For the converse, let $(\Om,R)$ be any open $\sigma$-trivial
regularization of the dynamical system $(A,\,G,\,\alpha)$.
Then $\Om=\Lambda\times_HG$, where $\Lambda$ is the image of
a section for $\Om/G$ and $H$ is the restriction of the stabilizer map to
$\Lambda$. Let $I=\ker R^{-1}(\Lambda)$ and let $B=A/I$. Then it is easily seen
that the canonical action of $\Om^S$ on $A$ given by Example \ref{ex1}
restricts to a subgroup action of $\Lambda^H$ on $B$, and it was shown in
\cite[Theorem 3]{ech2}
that there is a canonical $G$-equivariant isomorphism between
$A$ and $\Ind_{\Lambda^H}^GB$ which is given by
$$\(\Phi(a)\)(s)=\alpha_{s^{-1}}(a)+I.$$

One reason that induced systems are interesting in the investigation of crossed
products is the following result.

\begin{prop}[{\cite[Corollary 3]{ech2}}]
\label{propind}
Let $\beta$ be an action of $\Om^H$ on $B$. Then
$\Ind_{\Om^H}^GB\rtimes_{\Ind\beta}G$ is Morita equivalent to
$B\rtimes_{\beta}\Om^H$.
\end{prop}

%
%
\section{Twisted systems and induced
crossed products}
\label{sec:3}

Our main result in this section is the following theorem which
is a generalization of a result of Pedersen and Olesen
\cite{op} (which is essentially (1)~$\Longrightarrow$~(3) in the theorem
below  in case of ordinary Green twisting maps)
to generalized twisted systems.
In the following
we say that a covariant system $(A,\,G,\,\alpha)$ is
{\bem abelian}, if $G$ is an abelian group.

\begin{thm}\label{theotwist}
Suppose that $(A,\,G,\,\alpha)$ is an abelian dynamical system and
that $\Om$ is a locally compact Hausdorff space.
Let $N:\Om\to\frak K(G)$ be a continuous map.
Then the following conditions are equivalent.
\begin{enumerate}
\item There exists an open $G$-invariant
regularization $(\Om,R)$ for $(A,\,G,\,\alpha)$
for which there is a
\ggtm{}
with domain $\Om^N$.
\item There exists an open $\sigma$-trivial regularization
$\ghat{R}:\spec{A\rtimes_{\alpha}G}\to \Om\times_{N^{\perp}}\ghat{G}$
for the dual system $(A\rtimes_{\alpha}G,\,\ghat{G},\,\ahat{\alpha})$.
\item
The dual system $(A\rtimes_{\alpha}G,\,\ghat{G},\,\ahat{\alpha})$
is isomorphic to an induced system
$(\Ind_{\Om^{N^{\perp}}}^{\ghat{G}}B,\,\ghat{G},\,\Ind\beta)$,
for some \cs-algebra $B$ and some action $\beta$ of $\Om^{N^{\perp}}$
on $B$.
\end{enumerate}
Moreover, if the conditions \rom{(1)} to \rom{(3)} are satisfied, then
the algebra $B$ in \rom{(3)} may be chosen to be $A\rtimes_{\alpha,\dtau}G$,
where the action of $\beta^x$ of $N_x^{\perp}$ is given on the fiber
$A_x\rtimes_{\alpha^x,\dtau^x}G$ by the dual action of $N_x^{\perp}$
on this fiber.
\end{thm}

Note that the equivalence of (2) and (3) is a consequence of
\cite[Theorem 3]{ech2} (see also the remarks at the end of  \S1).
Thus in order to prove the equivalences of our theorem it is enough
to prove the implications (1)~$\Longrightarrow$~(2)
and (3)~$\Longrightarrow$~(1).
Before we start with the proofs we have to recall  some facts about
the space $\Rep(D)$ of equivalence classes of $*$-representations
(with bounded dimension) of a \cs-algebra $D$.

\begin{prop}\label{proprep}
Let $D$ be a \cs-algebra and let $(\pi_i)_{i\in I}$ be a net in $\Rep(D)$
which converges to $\pi\in\Rep(D)$.
\begin{enumerate}
\item If $\rho\in\Rep(D)$
is such that $\ker\rho\supseteq \ker\pi$, then
$\pi_i\to\rho$ in $\Rep(D)$.
\item Suppose that $\rho\in \hat{D}$
is such that $\ker\rho\supseteq\ker\pi$,
and let $\cal D_i\subseteq\hat{D}$ be
such that $\ker \cal D_i=\ker \pi_i$
for all $i\in I$. Then there
exists a subnet $(\pi_j)_{j\in J}$ and elements $\rho_j\in\cal D_j$
such that $\rho_j\to\rho$ in $\hat{D}$.
\end{enumerate}
\end{prop}
\begin{pf}
The first assertion is just \cite[Proposition 1.2]{fell6}, while the second
is a consequence of \cite[Theorem 2.2]{schoch} (note that Schochetman proves
his result only for separable $D$, but the same proof applies in the
general case).
\end{pf}

In the sequel we need an explicit description of $\cs(\Om^N)$
in case where $G$ is abelian.

\begin{lem}\label{COMN}
Suppose that $\Om$ is a locally compact space and $N:\Om\to\frak K(G);
x\mapsto N_x$
is a continuous map for some locally compact abelian group $G$.
Then the map
$$M: \Om\times\ghat{G}\to \specnp{\cs (\Om^N)}; (x,\chi)\mapsto (x,\chi|_N)$$
defines a homeomorphism between $\Om\times_{N^{\perp}}\ghat{G}$
and $\specnp{\cs (\Om^N)}$. As a consequence, $\cs (\Om^N)$ is isomorphic
to $C_0(\Om\times_{N^{\perp}}\ghat{G})$.
\end{lem}
\begin{pf}
It is clear that $M$ defines a bijection between
$\Om\times_{N^{\perp}}\ghat{G}$ and $\cs (\Om^N)$, so it is enough to
show that $M$ is continuous and open. The continuity is a direct
consequence of continuity of restricting representations (see
\cite[Proposition 7]{ech1}). For proving the openness, let
$\set{(x_i,\chi_i)}_{i\in I}$ be a net in $\Om\times G$ such that
$\set{(x_i,\chi_i|_{N_{_i}})}$ converges to $(x,\chi|_{N_x})$ in
$\specnp{\cs (\Om^N)}$.
Then $x_i\to x$ in $\Om$.
By the continuity of induction \cite[Proposition 6]{ech1} we
know that $\ind_{N_{x_i}}^G(\chi_i|_{N_{x_i}})$ converges to
$\ind_{N_x}^G(\chi|_{N_x})$ in $\Rep \(\cs (G)\)$. But it is well known
that $\ind_{N_{x_i}}^G(\chi_i|_{N_{x_i}})$ is weakly equivalent to the coset
$\chi_iN_{x_i}^{\perp}\subseteq\ghat{G}$ for each $i\in I$, and that
$\ind_{N_x}^G(\chi|_{N_x})$ is weakly equivalent to $\chi N_x^{\perp}$
(see, for example, \cite[Lemma~5.1]{dpw2}).
Thus it is a consequence of Proposition \ref{proprep} that we may pass
to a subnet in order to find elements $\mu_i\in N_{x_i}^{\perp}$ such that
$(x_i,\chi_i\mu_i)$ converges to $(x,\chi)$ in $\Om\times\ghat{G}$.
But this shows the openness of $M$.
\end{pf}

We use this lemma to prove

\begin{lem}\label{lem1}
Suppose that $(\Om,R)$ is a $G$-invariant open regularization
for the covariant system $(A,\,G,\,\alpha)$ and that $r:E\to\Om$ is the
corresponding \cs-bundle such that $A=\Gamma_0(E)$.
Assume that $\dtau$ is a \ggtm{} for $(A,\,G,\,\alpha)$
with domain $\Om^N$.
Then $A\rtimes_{\alpha}\Om^N$ is isomorphic to $A
\rtimes_{\id}\Om^N$,
the subgroup crossed product of $A$ with $\Om^N$ by the trivial action.
If $G$ is abelian, then this implies that $A\rtimes_{\alpha}\Om^N$
is canonically
isomorphic to the balanced tensor product
$A\otimes_{C_0(\Om)}C_0(\Om\times_{N^{\perp}}\ghat{G})$.
The corresponding homeomorphism of
$$\spec{A\otimes_{C_0(\Om)}C_0(\Om\times_{N^{\perp}}\ghat{G})}=
\set{\(\rho, (x,\chi)\); R(\rho)=x}\subseteq
\hat{A}\times (\Om\times_{N^{\perp}}\ghat{G})$$
onto $\spec{A\rtimes_\alpha\Om^N}$
is given by the map
\begin{equation}
\label{eq:specmap}
\(\rho, (x,\chi)\)\mapsto \(x,\rho\times \chi(\rho\circ \dtau^x)\),
\end{equation}
where the latter term denotes
the representation of $A\times_\alpha\Om^N$
given by
$\rho\times \chi(\rho\circ \dtau^x)$ acting
on the quotient $A_x\rtimes_{\alpha^x}N_x$.
\end{lem}
\begin{pf}
Let $q:\Om^N\to\Om$ denote the projection. Then
it follows by straight forward computations that the map
$$\Psi:\Gamma_c(q^*E)\to\Gamma_c(q^*E); \Psi(f)(x,n)=f(x,n)\dtau^x_n$$
extends to an isomorphism from $A\rtimes_{\alpha}\Om^N$
onto  $A\rtimes_{\id} \Om^N$.
Suppose now that $G$ is abelian.
Consider the map
$\Phi: A\odot C_c(\Om^N)\to A\rtimes_{\id}\Om^N$
which is defined by
$$\Bigl(\Phi\(\sum_{i=1}^na_i\otimes f_i\)\Bigr)(x,n)=
\sum_{i=1}^na_i(x)f_i(x,n).$$
It is easily seen that $\Phi$ extends to a surjective $*$-homomorphism
from $A\otimes \cs (\Om^N)$ onto $A\rtimes_{\id}\Om^N$ such that the
kernel of this map is exactly the intersection of
all representations $\(\rho, (x,\chi)\)\in \hat{A}\times\specnp{ \cs (\Om^N)}$
such that $R(\rho)=x$. This implies that $\Phi$ factors through an isomorphism
between $A\otimes_{C_0(\Om)}\cs (\Om^N)$ and $A\rtimes_{\id}\Om^N$.
Thus we conclude from Lemma \ref{COMN} that $A\rtimes_{\alpha}\Om^N$
is canonically isomorphic to
$A\otimes_{C_0(\Om)}C_0(\Om\times_{N^{\perp}}\ghat{G})$.

Finally, let $\(\rho, (x,\chi)\)\in
\spec{A\otimes_{C_0(\Om)}C_0(\Om\times_{N^{\perp}}\ghat{G})}$,
and let $a\in A$ and $f\in C_c(\Om^N)$, the latter viewed as an element of
$C_0(\Om\times_{N^{\perp}}\ghat{G})$ via Fourier transform.
Then we compute
\begin{gather*}
\(x, \rho\times \chi(\rho\circ \frak t^x)\)\(\Psi^{-1}\Phi(a\otimes f)\)
=\rho\times \chi(\rho\circ \frak t^x)((a(x)f(x,\cdot)(\frak t^x)^*)
\\
\qquad=\int_{N_x}\rho\(a(x)\)f(x,n)\chi(n)\,d_{N_x}n=
\(\rho, (x,\chi)\)(a\otimes f),
\end{gather*}
which finishes the proof.
\end{pf}

In the sequel, it will be necessary to use some continuity results for
inducing and restricting representations, which are not explicitly stated in
\cite{ech1}.

\begin{prop}\label{propindres}
Suppose that $(\Om,R)$ is a $G$-invariant open regularization of the
system $(A,\,G,\,\alpha)$
and let  $N:\Om\to\frak K(G)$ be a continuous map such that
$N_x$ is normal in $G$ for all $x\in \Om$. Then the maps
$$\ind:\cal S(A\rtimes_{\alpha}\Om^N)\to \Rep(A\rtimes_{\alpha}G);
(x,\rho\times V)\mapsto\ind_{N_x}^G(\rho\times V)$$
and
$$\res:\spec{A\rtimes_{\alpha}G}\to\cal S(A\rtimes_{\alpha}\Om^N);
\pi\times U\mapsto(P(\pi\times U), \pi\times U|_{N_{P(\pi\times U)}})$$
are continuous, where $P:\spec{A\rtimes_{\alpha}G}\to\Om$ denotes the
canonical projection.
\end{prop}
\begin{pf}
By Proposition \ref{propisom1} we may view $\cal S(A\rtimes_{\alpha}\Om^N)$
topologically as a subset of $\cal S(\Om^N,A,\alpha)$.
But it is a direct consequence of \cite[Corollary 2]{ech1} that the
extension of $\ind$ to $\cal S(\Om^N,A,\alpha)$ is continuous.
Thus our map is continuous, too. The continuity of $\res$
follows easily from Proposition \ref{propisom1} and \cite[Proposition 7]{ech1}
by identifying $A\rtimes_{\alpha}G$ with $A\rtimes_{\alpha}\Om^G$ as in
Lemma~\ref{lemid}.
\end{pf}

We now prove the openness of the canonical projection
from $\spec{A\rtimes_{\alpha,\dtau}G}$ onto $\Om$ when $\dtau$ has domain
$\Om^N$. In fact we will prove a slightly more general result which we
will need in the proof of \thmref{theotwist}.

\begin{prop}\label{proptwistreg}
Suppose that $(\Om, R)$ is a $G$-invariant open regularization,
$N:\Om\to\frak K(G)$ a continuous map, and
$\dtau$ a
\ggtm{}  for $(A,\,G,\,\alpha)$
with domain $\Om^N$. Suppose further that $(\Lambda,P)$ is another
$G$-invariant open regularization of $(A,\,G,\,\alpha)$ such that
there exists an open continuous map $q:\Lambda\to\Om$ satisfying
$R=q\circ P$. Let $Q:\spec{A\rtimes_{\alpha,\dtau}G}\to\Lambda$
denote the restriction  to $\spec{A\rtimes_{\alpha,\dtau}G}$
of the canonical projection from
$\spec{A\rtimes_{\alpha}G}$ onto $\Lambda$.
Then $Q$ is open and surjective whenever $G/N_x$ is amenable for all
$x\in\Om$.
\end{prop}
\begin{pf}
Recall that $Q(\pi\times U)=y$ if and only if $\ker\pi\supseteq P^{-1}(
\set{y})$.
This implies that
$\pi\times U$ factors through the quotient
$A_y\rtimes_{\alpha^y,\dtau^y}G$, where
$\alpha^y$ denotes the canonical action of $G$ on $A_y=A/\ker P^{-1}(
\set{y})$
and $\dtau^y$ denotes the image of $\dtau^{q(y)}$ under the canonical
quotient map $A_{q(y)}\to A_y$ (which exists due to the fact that
$\ker R^{-1}(
\set{q(y)})\subseteq P^{-1}(
\set{y})$). Now let $\rho$ be
a faithful representation of $A_y$. Then we may produce a
covariant representation $\rho\times V$ of
$A_y\rtimes_{\alpha^y,\dtau^y}N_{q(y)}$
by defining $V=\rho\circ \dtau^y$. Then $\rho\times V$ preserves the
twist $\dtau^y$, and
$\ind_{N_{q(y)}}^G(\rho\times V)$ is a representation of
$A_y\rtimes_{\alpha^y}G$ which preserves $\dtau^y$ by
\cite[Corollary 5]{green1}.
It follows from \cite[Proposition 13]{green1} and the amenability of
$G/N_{q(y)}$ that
this induced representation defines a faithful representation of
$A\rtimes_{\alpha^y,\dtau^y}G$. Regarding
$\ind_{N_{q(y)}}^G(\rho\times V)$ as a representation of
$A\rtimes_{\alpha,\dtau}G$, we have
$\ker\(\ind_{N_{q(y)}}^G(\rho\times V)\)=\ker Q^{-1}(
\set{y})$.

Now let $\set{y_i}_{i\in I}$ be a net in $\Lambda$ such that
$y_i\to y$  and let
$\pi\times U\in Q^{-1}(\set{y})$.
We have to show that there exists a subnet
$\set{y_j}_{j\in J}$ of $\set{y_i}_{i\in I}$, and a net
$\set{\pi_j\times U^j}_{j\in J}\subseteq
\spec{A\rtimes_{\alpha,\dtau}G}$,
converging to $\pi\times U$,
such that $Q(\pi_j\times U^j)=y_j$ for all $j\in J$.

For each $i\in I$, let $\rho_i\in \Rep(A)$ be such that
$\ker\rho_i=\ker P^{-1}(\set{y_i})$, and let $\rho\in\Rep(A)$
be such that
$\ker\rho=\ker R^{-1}(\set{y})$. Then the same arguments as used in the proof
of
\cite[Proposition 11]{ech2}
show that $\rho_i\to\rho$ in $\Rep(A)$.

Looking at the subgroup crossed product $A\rtimes_{\id}\Om^N$
by the trivial action of $\Om^N$ on $A$, we observe that
$\(q(y_i), \rho_i\times 1_{\cal H_{\rho_i}}\)$ converges to
$(q(y),\rho\times 1_{\cal H_{\rho}})$ in $\cal S(A\rtimes_{\id}\Om^N)$,
where $1_{\cal H_{\rho_i}}$ denotes the trivial representation of
$N_{q(y_i)}$ into $\cal H_{\rho_i}$ and $1_{\cal H_{\rho}}$ is defined
similarly. If $\Phi:A\rtimes_{\alpha}\Om^N\to A\rtimes_{\id}\Om^N$ is given as
in Lemma \ref{lem1}, then we see that the representations
$\(q(y_i), \rho_i\times 1_{\cal H_{\rho_i}}\)\circ \Phi$
of $A\rtimes_{\alpha}\Om^N$ are just given by
the subgroup representations $\(q(y_i),\rho_i\times V^i\)$
with $V^i=\rho_i\circ \dtau^{q(y_i)}$.

It follows that $(q(y_i), \rho_i\times V^i)$
converges to $(q(y),\rho\times V)$
in $\cal S(A\rtimes_{\alpha}\Om^N)$. Now the continuity of induction
implies that $\ind_{N_{q(y_i)}}^G(\rho_i\times V^i)$ converges to
$\ind_{N_{q(y)}}^G(\rho\times V)$ in $\Rep(A\rtimes_{\alpha}G)$.
But as above, we know that
$\ker\(\ind_{N_{q(y_i)}}^G(\rho_i\times V^i)\)=\ker Q^{-1}(\set{y_i})$ for each
$i\in I$. Thus by Proposition \ref{proprep} we may pass to
a subnet in order to find elements $\pi_i\times U^i\in Q^{-1}(\set{y_i})$
such that $\pi_i\times U^i$ converges to $\pi\times U$ in
$A\rtimes_{\alpha,\dtau}G$. This finishes the proof
\end{pf}

\begin{cor}\label{coropen}
Suppose that $\dtau$ is a \ggtm{} for
$(A,\,G,\,\alpha)$ with domain $\Om^N$ such that $G/N_x$ is amenable
for all $x\in\Om$. Then $A\rtimes_{\alpha,\dtau}G$ is isomorphic
to the section algebra $\Gamma_0(D)$ of a \cs-bundle
$q:D\to \Om$ with fibers $A_x\rtimes_{\alpha^x,\dtau^x}G$.
\end{cor}
\begin{pf} The Corollary follows from
Lee's theorem
\cite[Theorem~4]{lee2} and by taking $\Lambda=\Om$
in Proposition~\ref{proptwistreg}.
\end{pf}

Our next proposition gives the implication
$(1)\Longrightarrow(2)$
in Theorem \ref{theotwist}.

\begin{prop}\label{prop:1to2}
Suppose that $(\Om,R)$ is an open $G$-invariant regularization
for an abelian dynamical system $(A,\,G,\,\alpha)$, and that $N:\Om
\to\frak K(G)$ is a continuous map such that there is a
\ggtm{} $\dtau$ with domain $\Om^N$.
\begin{enumerate}
\item
If we identify the spectrum of $A\rtimes_\alpha\Om^N$ with $\set
{\(\rho,(x,\chi)\):R(\rho)=x}$ (via \lemref{lem1}), then $\tilde
R\(\rho,(x,\chi)\)=(x,\chi)$ defines an open $\ghat G$-equivariant
and $G$-invariant regularization of the
iterated system $(A\rtimes_\alpha\Om^N,\, G,\,\gamma^N)$
as defined in Example~\ref{ex3}.
\item
If $\pi\times U$ is an irreducible representation of
$A_x\rtimes_\alpha G$, then there is a unique $(x,\chi)\in\Om
\times_{N^\perp}\ghat{G}$ with $\ker(\pi\times U\restr{N_x})\supseteq
\ker\(\tilde R^{-1}\(\set{(x,\chi)}\)\)$.  In particular, $\ghat
R(\pi\times U)=(x,\chi)$ defines an open $\sigma$-trivial
regularization for the dual system $(A\rtimes_\alpha G, \,\ghat G,\,
\ahat \alpha)$.
\end{enumerate}
\end{prop}
\begin{pf}
Let $(\gamma^N,\dtau^N)$ denote the canonical generalized twisted
action of $G$ on $A\rtimes_{\alpha}\Om^N$ as given in Example \ref{ex3}.
We show  that there is a continuous open $\ghat{G}$-equivariant
and $G$-invariant surjection
$\tilde{R}$ from
$\spec{A\rtimes_{\alpha}\Om^{N}}$ onto
$\Om\times_{N^{\perp}}\ghat{G}$.
By Lemma \ref{lem1} we know that $\spec{A\rtimes_{\alpha}\Om^N}$
is homeomorphic to $\set{\(\rho,
(x,\chi)\); R(\rho)=x}\subseteq \hat{A}\times
(\Om\times_{N^{\perp}}\ghat{G})$ via the homeomorphism
$\(\rho, (x,\chi)\)\mapsto \(x,\rho\times \chi (\rho\circ
\dtau^x)\)$.
We define $\tilde{R}\((x,\rho\times \chi (\rho\circ \dtau^x)
)\)= (x,\chi)$.
It follows directly from the openness and surjectivity of $R$ that
$\tilde{R}$ is open and surjective, too.
It is also clear that $\tilde{R}$ is $\ghat{G}$-equivariant,
since on both spaces
the action of $\ghat{G}$ is given by pointwise multiplication of
characters.
In order to show that $\tilde{R}$ is $G$-invariant note first that
by the definition of a Green twisting map
$\alpha^x_s\circ \dtau^x=\dtau^x$
for all $x\in \Om$ and $s\in G$. Thus we conclude for $f\in \Gamma_c(q^*E)$:
\begin{multline*}
\(x, \rho\times \chi(\rho\circ \dtau^x)\)\(\gamma^N_s(f)\)=
\int_{N_x}\rho\(\alpha_s\(f(x,n)\)\dtau^x_n\)\chi(n)\,d_{N_x}n
\\
=\int_{N_x}\rho\(\alpha_s\(f(x,n)\dtau^x_n\)\)\chi(n)\,d_{N_x}n
=\(x, \rho\circ\alpha_s\times \chi(\rho\circ\alpha_s\circ \dtau^x)\)(f).
\end{multline*}
Thus we see that $\(x, \rho\times \chi(\rho\circ \dtau^x)\)\circ\gamma_s$
is equal to
$\(x, \rho\circ\alpha_s\times \chi(\rho\circ\alpha_s\circ \dtau^x)\)$
which has the same image under $\tilde{R}$.

We have seen above that
$\tilde{R}:\spec{A\rtimes_{\alpha}\Om^N}\to\Om\times_{N^{\perp}}\ghat{G}$
is an open $G$-invariant regularization for the system
$(A\rtimes_{\alpha}\Om^N,\,G,\,\gamma^N)$ such that the projection, say $r$,
from $\Om\times_{N^{\perp}}\ghat{G}$ onto $\Om$  satisfies the
relation $r\circ \tilde{R}=P$, where $P$ denotes the projection
from $\spec{A\rtimes_{\alpha}\Om^N}$ onto $\Om$.
Thus it follows from Proposition \ref{proptwistreg} that there is a canonical
open map $Q$ from $\spec{(A\rtimes_{\alpha}\Om^N)\rtimes_{\gamma^N,\dtau^N}G}$
onto $\Om\times_{N^{\perp}}\ghat{G}$,
and we define $\ghat{R}$ as the corresponding
open map from $\spec{A\rtimes_{\alpha}G}$ onto $\Om\times_{N^{\perp}}\ghat{G}$
via the isomorphism $\Psi$ between  $A\rtimes_{\alpha}G$
and $(A\rtimes_{\alpha}\Om^N)\rtimes_{\gamma^N,\dtau^N}G$ constructed
in Example \ref{ex3}.

Now let $\pi\times U\in \spec{A\rtimes_{\alpha}G}$. Then $\pi\times U$
factors through some quotient $A_x\rtimes_{\alpha^x}G$,
and is transported to the representation
$(\pi\times U|_{N_x})\times U$ of the fiber
$(A_x\rtimes_{\alpha^x}N_x)\rtimes_{\gamma^{N_x},\dtau^{N_x}}G$
of $(A\rtimes_{\alpha}\Om^N)\rtimes_{\gamma^N,\dtau^N}G$ via
$\Psi$.
Thus we see that $\ghat{R}(\pi\times U)=(x,\mu)$ for some
 $(x,\mu)\in \Om\times_{N^{\perp}}\ghat{G}$
 if and only if the kernel of the representation
$(x,\pi\times U|_{N_x})$ of $A_x\rtimes_{\alpha^x}N_x$ contains
$\ker \tilde{R}^{-1}(\set{(x,\mu)})$.
Now the dual action of $\chi\in\ghat{G}$ on the element $\pi\times U
\in \spec{A\rtimes_{\alpha}G}$ is given by multiplication of $\chi$ with $U$,
i.e., $(\pi\times U)\circ\ahat{\alpha}_{\chi^{-1}}=\pi\times\chi U$.
This representation corresponds to
the representation $(\pi\times\chi U|_{N_x})\times\chi U$
of the fiber $(A_x\rtimes_{\alpha^x}N_x)\rtimes_{\gamma^{N_x},\dtau^{N_x}}G$.
Hence it follows from the $\ghat{G}$-equivariance of $\tilde{R}$
that $\ghat{R}(\pi\times\chi U) =(x,\chi\mu)$ if and only if
$Q(\pi\times U)=(x,\mu)$.
Thus we see that $(\Om\times_{N^{\perp}}\ghat{G},\ghat{R})$ is a
$\sigma$-trivial open regularization for the dual system
$(A\rtimes_{\alpha}G,\,\ghat{G},\,\ahat{\alpha})$.
\end{pf}

Since $G$ is abelian,
the continuity of the map $N:\L\to\frak K(G)$ implies that
$N^\perp
:\L\to\frak K(\ghat G)$,
defined by $N^\perp_x=(N_x)^\perp$, is
also continuous \cite{dpw4}.
The implication $(3)\Longrightarrow(1)$ in Theorem~\ref{theotwist}
will follow from the
next result.

\begin{prop}\label{threeone}
Suppose that $(A,\,G,\,\alpha)$ is an abelian dynamical system.
Assume that $q:E'\to\L$ is a \cs-bundle such that $D=\Gamma_0
(E')$ admits an action $\gamma$ of $\L^{N^\perp}$ on $D$ such
that there is a covariant isomorphism from
$(A\rtimes_\alpha G,\,\ghat G,\,\ahat\alpha)$ onto
$(\ind_ {\L^{N^\perp}}^{\ghat G}D,\,\ghat G,\,\ind\gamma)$.  Then there is
an open $G$-invariant regularization for
$(A,\,G,\,\alpha)$ for which there is a \ggtm{} with domain $\L^N$.
\end{prop}

The proof of the above proposition will depend heavily on the
following two lemmas which may be of independent interest.

\begin{lem}\label{lemthreeone}
Suppose that $\alpha=\set{\alpha^x}_{x\in\L}$ is an action of
$\L^N$ on a \cs-algebra $D$, and that $(A,\,G,\,\beta)$ is the
induced system $(\ind_{\L^N}^G D,\,G,\,\ind\alpha)$.
Then there is a $\ghat G$-invariant open
regularization $(\ghat P,\L)$ of the dual system
$(A\rtimes_\beta G,\,\ghat G,\,\ahat \beta)$
with respect to which there is a
\ggtm{} for $(A\rtimes_\beta
G,\,\ghat G,\,\ahat\beta)$ with domain $\L^{N^\perp}$.
\end{lem}

\begin{pf}
It follows from \cite[Theorem~3]{ech2} that there is an open
$G$-invariant regularization $(P,\L)$ for
$(A,\,G,\,\beta)$.  As a consequence of
\cite[Theorem~2.1]{dpw5}, there is an open
surjection $\ghat P:
\Prim(A\rtimes_\beta G)\to\L$, and this map is easily seen to be
$\ghat G$-invariant.

In particular,
we may assume that $A=\G_0(E')$ and that $D=\G_0(E)$ where $E$ and
$E'$ are \cs-bundles over $\L$.  Furthermore, $A_x=\ind
_{N_x}^G D_x$.  Thus if $\sigma\in N_x^\perp$, then we can
define $\tilde \dtau_\sigma^x\in\UM\(A_x\)$ by
$$(\tilde \dtau_\sigma^x f)(s)=\overline{\sigma(s)}f(s) = (f\tilde
\dtau_\sigma^x)(s).$$
We let $\dtau_\sigma^\l$ be the image of $\tilde \dtau_\sigma^\l$ under
the natural inclusion of $\UM\(A_x\)$ into $\UM(A_x\rtimes_{
\alpha^\l}G)$.  Thus, if we view elements of $C_c(G,A)$ as
functions on $G\times\L\times G$, we have
\begin{equation*}
\(\dtau^\l_\sigma f(s)(\l)\)(t) = \overline{\sigma(t)}f(s)(\l)(t)
\quad\text{and}\quad
(f(s)(\l)\dtau_\sigma^\l)(t)=
\overline{\sigma(s^{-1}t)}f(s)(\l)(t).
\end{equation*}
The first step will be to show that $(\l,\sigma)\mapsto
\dtau_\sigma^\l f(\cdot)(\l)$ and $(\l,\sigma)\mapsto
f(\cdot)(\l)\dtau_\sigma^\l$ are continuous from $\L^{N^\perp}$ into
$E$ for a suitable collection of $f\in C_c(G,A)$.  Certainly we
may assume that $f(s)(\l)(t)=g(s)a(t)(\l)$ for $g\in C_c(G)$ and
$a\in A$.  So it suffices to see that $(\l,\sigma)\mapsto \tilde
\dtau_\sigma^\l a(\l)$ is continuous.  Using the definition of the
topology on $E$ (see \cite[Proposition 1.6]{fell4}),
 it will suffice
to show that given any $(\l_0,\sigma_0)\in\L^{N^\perp}$ and some
$b\in A$ satisfying $b(\cdot)(\l_0)=\tilde
\dtau_{\sigma_0}^{\l_0}a(\cdot)(\l_0)$,
$$\| \tilde \dtau_\sigma^\l a(\cdot)(\l)-b(\cdot)(\l)\|_{A_x}$$
goes to zero as $(\l,\sigma)$ approaches $(\l_0,\sigma_0)$.
However,
$$\| \tilde \dtau_\sigma^\l a(\cdot)(\l)-b(\cdot)(\l)\|_{A_x}=
\sup_{t\in G}\|\overline{\sigma(t)}a(t)(\l)-b(t)(\l)\|_{D_x}.$$
Furthermore the map $(\l,\sigma,t)\mapsto\|\overline{ \sigma(t)}
a(t)(\l)-b(t)(\l)\|_{D_x}$
is continuous from $\L^{N^\perp}\times G$ into $\RR^+$.  (Just
consider the map from $\L\times\ghat G\times G$ with any $a,b\in
C_b(G,D)$.)

We may assume that $a$, and hence $b$, have compact support in
$\L\times_N G$.  Let $\kappa:\L\times G\to\L\times_N G$ be the
quotient map.  Observe that there is a compact set $C\subseteq
\L\times G$ such that $\kappa(C)\supseteq \supp a\cup \supp b$
\cite[Lemma~2.3]{dpw3}.
Now choose $\phi\in C_c(\L\times G)$ such that $0\le \phi \le 1$
and $\phi\equiv 1$ on $C$.  Then define $\psi(\l,\sigma,t)$ to be
$\phi(\l,t)\|\overline{\sigma(t)}a(t)(\l)-b(t)(\l)\|_{D_x}$.
The point being that $\|\overline{\sigma(t)}a(t)(\l)-b(t)(\l)\|_{D_x}$
depends only on the class of $t$ in $G/N_x$ so that
\begin{align*}
\|\tilde \dtau_\sigma^\l a(\cdot)(\l)-a(\cdot)(\l)\|_{A_x} &=
\sup_{t\in G}
\|\tilde \dtau_\sigma^\l a(t)(\l)-a(t)(\l)\|_{D_x} \\
&=\sup_{t\in G}\psi(\l,\sigma,t).
\end{align*}
But we may assume that $\sigma$ lies in a compact neighborhood of
$\sigma_0$, and, hence, that $\psi$ has compact support in $\L
\times \ghat G\times G$.
Then, if $\psi_{(\l,\sigma)}\in C_c(G)$ is defined by
$\psi_{(\l,\sigma)}(s)=\psi(\l,\sigma,t)$, the map $(\l,\sigma)
\mapsto \psi_{(\l,\sigma)}$ is continuous from $\L\times\ghat G$
to $C_0(G)$.  Since $\psi_{(\l_0,\sigma_0)}=0$, we have
$\|\psi_{(\l,\sigma)}\|_\infty<\epsilon$ for $(\l,\sigma)$ near
$(\l_0,\sigma_0)$ which is what we wanted.

It is now immediate that $\sigma\mapsto \dtau^\l_\sigma$ is strictly
continuous on $N_\l^\perp$, and it follows from straightforward
calculations that $\dtau^\l_{\sigma\rho}=\dtau^\l_\sigma \dtau^\l_\rho$
if $\sigma,\rho\in N_\l^\perp$, that
if $\sigma\in N_\l^\perp$, then $\Ad \dtau_\sigma^\l = \ahat\alpha^\l
_\sigma$, and that $\ahat\alpha^\l_\rho(\dtau^\l_\sigma)= \dtau
^\l_\sigma$ for all $\rho\in\ghat G$ and $\sigma\in N_\l^\perp$.
Thus $\sigma\mapsto \dtau^\l_\sigma$ is an ordinary Green twisting
map for $(A_x,\,\ghat G,\,\ahat \alpha^\l)$.
This completes the proof of the
lemma.
\end{pf}

\begin{lem}\label{lemmthreeone}
Suppose that $(A,\,G,\,\alpha)$ and $(B,\,G,\,\beta)$ are Morita
equivalent via $(X,u)$.  Let $h_X:\Prim(B)\to\Prim(A)$ be the
Rieffel homeomorphism.  If $(P,\L)$ is an open $G$-invariant
regularization of
$(B,G,\beta)$  such that there
is a
\ggtm{}
with domain $\L^N$, then $(P\circ h_X,\L)$ is an open
regularization of $(A,\,G,\,\alpha)$ for which there exists a
\ggtm{} with domain $\L^N$.
\end{lem}

\begin{pf}
Let $\dtau$ be a \ggtm{} for $\beta$.
We may assume that $B=\G_0(E)$
for a \cs-bundle $p:E\to\L$ induced by the map $P$ and that
$A=\G_0(E')$ for a \cs-bundle $E'$ induced by the map $P\circ h_X$.
Then $(A_x,\,G,\,\alpha^\l)$ is Morita equivalent to $(B_x,\,G,\,\beta^\l)$
via $(X^x,u^\l)$, where $X^x=X/X\cdot{\ker (P^{-1}(\set x))}$, and
$n\mapsto\dtau^\l_n$ is, by
definition, a Green twisting map for $(B_x,\,G,\,\beta^\l)$.  Then
\cite[Proposition~2]{ech5} implies that $\ds^\l$, defined by $
\ds^\l_n\cdot\xi=u^\l_n(\xi)\cdot\dtau^\l_n$, $\xi\in X^x$ defines a
Green twisting map for $\alpha^\l$.  Naturally, we claim that
$\ds(\l,n)=\ds^\l_n$ defines a \ggtm{}
for $(A,\,G,\,\alpha)$.  At this point, we need
only show that the maps $(\l,n)\mapsto \ds^\l_n a(\l)$
and $(\l,n)\mapsto a(\l)\ds^\l_n$ are continuous from $\L^N$ to
$E'$ for each $a \in A$.

Using \cite[Propositions 1.6~and 1.7]{fell4}, realize $X$ as the section
space of a Banach (space) bundle $q:F\to \L$.  We claim that the
maps from $F*F$ to $A$ and $B$, respectively, defined by $(a,b)
\mapsto \lip {A_{q(a)}}<a,b>$ and $(a,b)\mapsto \rip {B_{q(b)}} <a,b>$
are continuous.  To see this, consider a convergent net $(a_i,b_i) \to
(a,b)$ in $F*F$.  Suppose $q(a_i)=q(b_i)=\l_i$.
Choose sections $\xi,\eta\in X$ such that $\xi(\l)=a$ and $\eta(\l)=b$
(where, of course, $q(a)=\l=q(b)$).  Then $\|a_i-\xi(\l_i)
\|\to 0$ as does $\|b_i-\eta(\l_i)\|$.  Thus
$$\lip {A_{\l_i}}\,<a_i,b_i>\to\lip {A_{\l}}<\xi(\l),\eta(\l) >
=\lip {A_{\l}}<a,b>.$$
(And similarly for $\rip {B_x}<\cdot,\cdot>$.)

The hypotheses that $(\l,n)\mapsto \dtau^\l_n\cdot a(\l)$ and
$(\l,n)\mapsto a(\l)\cdot\dtau^\l_n$  are continuous into $E$
implies that $(\l,n)\mapsto \xi(\l)\cdot\dtau^\l_n$ is continuous
into $F$.  (Simply expand $\|\xi(\l)\cdot\dtau^\l_n-\eta(\l)\|^2$ in
terms of $\rip {B_x}<\cdot,\cdot>$.)

Finally, we may check continuity for $a\in A$ of the form $
a(\l)=\lip {A_x}<\xi(\l),\eta(\l)>$.  But then
$$\ds^\l_n\cdot a(\l)
=\lip {A_x}<\ds^\l_n\cdot \xi(\l), \eta(\l)>
= \lip {A_x}<u^\l_n\(\xi(\l)\)\cdot\dtau^\l_n,\eta(\l)>,$$
which is continuous in view of the above discussion.  A similar
argument applies to $a(\l)\ds^\l_n$ and completes the
proof of the lemma.
\end{pf}
\begin{remark}
Suppose that $X$ is an $A$--$B$-equivalence bimodule and that $A$
and $B$ are both section algebras of \cs-bundles over $\L$, say
determined by continuous open maps $P:\Prim(B)\to\L$ and $P':\Prim
(A)\to \L$.  There is, unfortunately, no reason to suspect
that
$P'=P\circ h_X$
in general.  However, note that the actions of $A$ and $B$ on $X$
define left and right $C_0(\L)$ module actions on $X$ (see
\cite[p.~187]{ra1}).  Then
$P'=P\circ h_X$ exactly when we have $\phi\cdot x=x\cdot\phi$ for all $x\in X$
and
$\phi\in C_0(\L)$. \end{remark}

\begin{pf*}{Proof of Proposition \ref{threeone}}
More or less by assumption, there is
a covariant $C_0(\L)$-isomorphism of the double-dual system
$((A\rtimes_\alpha G)\rtimes_{\ahat \alpha}\ghat
G,\,G,\,\aHat{\aHat\alpha})$ onto $(\indlgh\rtimes_{\ind\gamma}\ghat G, \,G,
\,\spec{\ind\gamma})$.  However, it follows from
\cite[Corollary~2]{ech5} that the
double-dual system is Morita equivalent to $(A,\,G,\,\alpha)$.
Therefore, $(\indlgh \rtimes_{\ind\gamma}
\ghat G,\,G,\,\spec{\ind\gamma})$ is Morita equivalent to $(A,\,G,\,\alpha)$.
The proposition now follows from Lemmas \ref{lemthreeone}~and
\ref{lemmthreeone}.
\end{pf*}
\begin{pf*}
{Proof of \thmref{theotwist}}
It is a consequence of Propositions \ref{prop:1to2}~and \ref{threeone}
as well as \cite[Theorem~3]{ech2} that conditions (1),~(2), and
(3) are equivalent.  Thus all that remains is the final
assertion.

So suppose that $\dtau$ is a
\ggtm{} for
$(A,\,G,\,\alpha)$ with domain $\Om^N$, and that $\tilde R$ and $\ghat
R$ are as in \propref{prop:1to2}.  Then \cite[Theorem~3]{ech2}
implies that the algebra $B$ in (3) is simply $(A\rtimes_\alpha G)/I$,
where $I=\ker\(\ghat R^{-1}\(\set{(x,1)\in\Om \times_{N^\perp}
G}\)\)$.  A moments reflection shows that it will suffice to show
that $I=I^\dtau$.

If $(x,\pi\times U)$ is an irreducible representation of
$A_x\rtimes_\alpha G$, then
\begin{equation*}
I=\bigcap\set{\ker(x,\pi\times U):\text{$x\in\Om$ and $
\ker(x,\pi\times U\restr{N_x})\supseteq \ker\(\tilde R^{-1}\(\set
{(x,1)}\)\)$}}.
\end{equation*}
But $\pi\times U$ preserves $\dtau^x$ if and only if $\pi\times
U\restr{N_x}$ does.  Although $\pi\times
U\restr{N_x}$ may not be irreducible, $\pi\times
U\restr{N_x}$ preserves $\dtau^x$ if and only if every
irreducible representation $\rho\times V$ of $A_x\rtimes_\alpha
N_x$ which satisfies $\ker(\rho\times V)\supseteq \ker(\pi \times
U\restr{N_x})$ preserves $\dtau^x$.
But if $\ker(\rho\times V)\supseteq\ker\(\pi\times U\restr{N_x})$,
then $\tilde R(\rho\times V)=(x,1)$.  Then \lemref{lem1} implies
that $V=\rho\circ\dtau^x$; that is, $\rho\times V$ preserves
$\dtau^x$.
It follows that $I=I^{\dtau}$ as claimed.
\end{pf*}

%
%
\section{Main Theorem}
\label{sec:4}

In this section $(A,\,G,\,\alpha)$ will be a {\it separable\/}
abelian dynamical system.  One of the basic invariants
employed in the study of crossed products is the Connes
spectrum $\Gamma(\alpha)$.  In our situation, Gootman has
shown \cite[Lemma~4.1]{goot2} that
\begin{equation*}
\Gamma(\alpha)=\bigcap \set{\ghat S_P:\text{$P$ is separated
in $\Prim(A\rtimes_{\alpha}G)$}},
\end{equation*}
where $\ghat S_P$ denotes the stabilizer of $P$ in $\ghat G$.
(Recall that a point $p$ in a topological space $X$ is called
separated if for all $q\in X\setminus\overline{\set{p}}$, $p$
and $q$ have disjoint \nbhd s.)  Actually,
we will be most interested in the subgroup $\ttg(\alpha)$ as
defined in \cite[\S5]{op}.  Note that
\begin{equation*}
\ttg(\alpha) = \bigcap \set{\ghat S_P : P\in
\Prim(A\rtimes_{\alpha}G)}.
 \end{equation*}
It is the fact that one can have
$\ttg(\alpha) \subsetneqq \Gamma(\alpha)$ even
for $G$-simple systems which makes characterizing simple
crossed products so subtle (see, for example,
 \cite[Theorems 5.5~\& 5.7]{op});
it also makes it necessary to employ $\ttg$
rather than $\Gamma$ in our next theorem (but see
Remark~\ref{rem:connes}).
\begin{thm}\label{thm:strivial}
Suppose that $(A,\,G,\,\alpha)$ is a separable abelian dynamical
system.
Then, with respect to the dual $\ghat G$-action, $\Prim(\acg)$
is a $\sigma$-trivial $\ghat G$-space if and only if the
following three conditions are satisfied.
\begin{enumerate}
\item
The quasi-orbit space $\Om=\qob G[\Prim(A)]$ is Hausdorff.
\item
If $\alpha^\w$ is the action on the quotient $A_\w$
corresponding to $\w\in\Om$ and if
$C_x=\ttg(\alpha^\w)^\perp$, then $C:\Om\to {\frak
K}(G)$ is continuous.
\item
There is a \ggtm{} for $(A,\,G,\,\alpha)$ with
domain $\Om^C$.
\end{enumerate}
 In this event, $\Prim(\acg)$ is $\ghat G$-homeomorphic to $\Om
 \times_{C^\perp}G$.
 \end{thm}
\begin{pf*}{Proof of necessity}
In the case that
$\Prim(\acg)$ is a $\sigma$-trivial $\ghat G$-space, then
certainly $\qob \ghat G [\Prim(\acg)] = \Prim(\acg)/\ghat G$ is
Hausdorff.  But $\Om=\qob G[\Prim(A)]$ is homeomorphic to
$\qob \ghat G [\Prim(\acg)]$ by \cite[Corollary~2.5]{gl}.
This establishes the necessity of~(1).

Now fix $\w\in\Om$, and consider the system
$(A_\w,G,\alpha^\w)$.  Since $\qob G[\Prim(A_\w)]
=\set{\text{pt}}$,
we employ \cite{gl} again to conclude that
$\qob \ghat G [\Prim(A_\w\rtimes_{\alpha^\w}G)] =
\set{\text{pt}}$.
Since $\Prim(A_\w\rtimes_{\alpha^\w}G) $ is a closed irreducible
subset of a $\sigma$-trivial space, it must
consist of a single orbit.
Therefore there is a single stabilizer group $\ghat S_\w=\ghat S_P$
for all $P\in\Prim(A_\w\rtimes_{\alpha^\w}G)$.  Since $P \mapsto
\ghat S_P$ is continuous by assumption, so is $\w\mapsto \ghat
S_\w^\perp=C_\w$ \cite{dpw4}.  Condition~(2) follows.

Using (1) and (2) above, we now have a continuous $\ghat
G$-ho\-me\-o\-mor\-phism $R:\Prim(\acg)\to \Om\times_{C^\perp}\ghat G$.
Now the necessity of~(3) is a corollary of $(2)\Longrightarrow(1)$
in \thmref{theotwist}.
\end{pf*}
\begin{remark}\label{rem:connes}
In the proof of necessity, we actually showed that
$C_\w=\Gamma(\alpha^\w)^{\perp}$
rather than $\ttg(\alpha^\w)^{\perp}$.  This follows because $\Prim(\acg)$,
and therefore $\Prim({A_\w\rtimes_{\alpha^\w}G})$, is Hausdorff;
hence \cite[Lemma~4.1]{goot2} implies that $\Gamma(\alpha^\w)=
\bigcap_{P\in\Prim(A_\w\rtimes_{\alpha^\w}G)}\ghat S_P$ in this case.
\end{remark}
\begin{pf*} {Proof of sufficiency}
Let $q:\Prim(A)\to \qob G[\Prim(A)]$ be the quasi-orbit map.
Recall that $q$ is continuous and open by
\cite[p. 221]{green1}.
Therefore, $(q,\Om)$
is a $G$-invariant regularization of $(A,\,G,\,\alpha)$, and by
assumption, there is
a \ggtm{} $\dtau$ with domain $\Om^C$ for $(A,\,G,\,\alpha)$.  By
\propref{prop:1to2} (or $(1)\Longrightarrow(2)$ of \thmref{theotwist}),
there is an open $\sigma$-trivial regularization $\ghat
R:\Prim(\acg) \to \Om\times_{C^\perp}\ghat G$ such that
the composition of
$R$ with the projection ${\operatorname{pr}}:\Om\times_{C^{\perp}}\ghat
G\to\Om$
is the open surjection $P:\Prim(\acg)\to\Om$
arising from $q$ via \cite{dpw5} and \cite{lee2}.

It will suffice to prove that $\ghat R$ is a homeomorphism---which
amounts to showing that $\ghat R$ is injective.  Since
$$
\begin{diagram}
\node{\Prim(A)} \arrow[2]{e,t}{\ghat R} \arrow{se,b}{P}
\node[2]{\Omega\times_{C^\perp}\ghat G}
\arrow{sw,b}{\operatorname{pr}} \\
\node[2]{\Om}
\end{diagram}
$$
commutes, it will suffice to consider a single fiber.  Thus we
fix $\w\in\Om$ and consider $\ghat R_\w:\Prim(A_\w\rtimes_{\alpha^\w}
G)\to \ghat G/C_\w^\perp\cong \ghat C_\w$.
Since $\ghat R_\w$ is $\ghat G$-equivariant and surjective, it will
suffice to show that $\Prim(A_\w\rtimes_{\alpha^\w}G)$ is
homeomorphic to $\ghat C_\w$ as a $\ghat G$-space.

Since $\dtau^\w$ is an (ordinary) Green twisting map for $(A_\w,\,G,\,
\alpha^\w)$ with domain $C_x^\perp$, we have from \cite{op}  or
\thmref{theotwist} that $A_\w\rtimes_{\alpha^\w}G$ is covariantly
isomorphic to $\Ind_{C^\perp}^{\ghat G}(A_\w\rtimes_{\alpha^\w,\dtau^\w}
G)$, where the action on $A_\w\rtimes_{\alpha^\w,\dtau^\w}G$ is
the dual action.  But $A_\w\rtimes_{\alpha^\w,\dtau^\w}G$ is
simple by \cite[Theorem~5.7]{op}.  An application of
\cite[Proposition~10]{ech2} completes the proof.
\end{pf*}

Recall that $(A,\,G,\,\alpha)$ is called {\it regular\/} \cite[p.~223]{green1}
if each quasi-orbit $Q\in\qob G[\Prim(A)]$ is locally closed, and
there is a $P\in Q$ such that
the map $R_P$ defined by $s\cdot S_P\mapsto \alpha_s(P)$
defines a homeomorphism of $G/S_P$ onto $Q$,
$S_P$ denoting the stabilizer of $P$.
(This implies that
the quasi-orbit space coincides with the orbit space.)
Also recall that if $A$ is type I and $\pi\in\hat A$,
then the symmetrizer subgroup
$\Sigma_{\ker\pi}=\Sigma_{\pi}$ \cite{go,horr,goot2}
is the image in $G$ of the center of
the extension of $S_\pi=S_{\ker\pi}$
by $\T$ determined by the Mackey obstruction
at $\pi$.
%
\begin{cor}\label{cor:symm}
If $A$ is type~I and $(A,\,G,\,\alpha)$ is regular, then
condition~{\normalshape (2)} in \thmref{thm:strivial} can be replaced by
\begin{description}
\item[$(2)'$]
The symmetrizer map $P\mapsto \Sigma_P$ is continuous, constant
on quasi-orbits, and defines a continuous map $C:\Om\to\frak
K(G)$.
\end{description}
\end{cor}
\begin{pf}
Since $A$ is type~I, the symmetrizer map will be constant on
quasi-orbits by \cite[Theorem~2.3]{goot2}. Hence we obtain a
well defined map
$\Sigma:\Om\to \frak K(G)$.
Thus we only have to
check that $\ghat S_P=\Sigma_{\ghat R(P)}^\perp$ for each
$P\in\Prim(\acg)$,
where $\ghat R:\Prim(\acg)\to \Om$ denotes the canonical map.
By the regularity assumption we know that
$\Prim(A_x)$ is homeomorphic to $G/S_J$ via
$sS_J\mapsto \alpha_s(J)$,
where $J$ denotes any element in the quasi-orbit $x$.

So let $J\in x=\ghat{R}(P)$.
Then \cite[Theorem~17]{green1} implies that $A_\w \rtimes_{
\alpha^\w} G$ is (strongly) Morita equivalent to $(A_\w/J)\rtimes
_{\alpha^\w}S_J$.  Thus a typical primitive ideal in $A_x\rtimes_{\alpha^x}G$
is the kernel of a representation of the form $\Ind_{S_J}^G(\pi
\times U)$ where $\pi\times U$ is an irreducible representation
of $(A_\w/J)\rtimes_{\alpha^\w}S_J$.  Furthermore, if $\ahat
\alpha^\w$ is the dual action of $\ghat G$ on $A_\w\rtimes_{\alpha^\w}
G$ and $\beta$ is the dual action of $\spec{S_J}$ on
$(A_\w/J)\rtimes_{\alpha^\w}S_J$, then \cite[Lemma~2.3]{rr}
implies that
\begin{equation*}
\Ind_{S_J}^G(\pi\times U)\circ({\ahat\alpha}^\w_\gamma)^{-1}=
\Ind_{S_J}^G(\pi\times\gamma\restr{S_J}U) =
\Ind_{S_J}^G\((\pi\times U)\circ\beta_{\gamma\restr{S_J}}^{-1}\)
\end{equation*}
It follows that the isotropy group $\ghat S_P$ for the dual action
at $P=\ker\(\Ind_{S_J}^G(\pi\times U)\)$ is
\begin{equation}
\label{eq:stabgp}
\set{\gamma\in \ghat G: \gamma\restr{S_J}\in \ghat S_{\ker(\pi\times
U)}},
\end{equation}
where $\ghat S_{\ker(\pi\times
U)}$ is the stabilizer group for the dual action of $\spec{S_J}$ on
$(A_\w/J)\rtimes_{\alpha^\w} S_J$ at $\ker(\pi\times U)$.  But,
since $A$ is type~I, $(A_\w/J)\cong\K$ and it is well known that
in this situation $\Prim\((A_x/J)\rtimes_{\alpha^x}S_J\)$ is
$\spec{S_J}$-homeomorphic to $\spec{\Sigma_J}$ (see for instance
\cite{horr}). This implies that
$\ghat S_P=\Sigma_J^{\perp}=\Sigma_{R(P)}^\perp$
  which completes the proof.
\end{pf}

\begin{remark}
In a forthcoming paper \cite{ech4} the first author shows that for
any separable abelian system $(A,\,G,\,\alpha)$ such that $A$ and
$A\rtimes_{\alpha}G$ are type I the $\ghat{G}$-stabilizer of a primitive
ideal $P$ of $A\rtimes_{\alpha}G$ is always equal to
the annihilator of the symmetrizer $\Sigma_{\Res P}$,
where $\Res :\Prim(A\rtimes_{\alpha}G)\to\cal Q^G\(\Prim(A)\)$
is the canonical restriction map as given by \cite[Corollary 19]{green1}.
Using this we see that Corollary \ref{cor:symm} holds under
the assumption that $A$ and $A\rtimes_{\alpha}G$ are type I
without assuming that $(A,G,\alpha)$ is regular.
\end{remark}

A rather special case of regular actions, but one which is still
of considerable interest, is the case where the induced $G$-action
is trivial on the spectrum.  In the type~I case, our theorem
specializes as follows.
\begin{cor}\label{cor:unitary}
Suppose that $(A,\,G,\,\alpha)$ is a separable abelian dynamical
system such that $A$ is type~I and the action of $G$ on $\Prim(A)$
is trivial.  Then $\Prim(\acg)$ is a $\sigma$-trivial $\ghat G$-space
if and only if the following three conditions are satisfied.
\begin{enumerate}
\item
$\Om=\Prim(A)$ is Hausdorff.
\item
The symmetrizer map $\Sigma:\Om\to\frak K(G); x\mapsto \Sigma_x$ is continuous.
\item
The restriction of $\alpha$ to $\Om^\Sigma$ defines a unitary subgroup
action of $\Om^\Sigma$ on $A$.
\end{enumerate}
In particular, if $\hat A$ is Hausdorff, then $\spec{\acg}\cong_{\ghat
G}\Om
\times\ghat G$ if and only if $\alpha $ is unitary.
\end{cor}
\begin{pf}
Since a \ggtm{} is unitary on its domain (see Remark~\ref{rem:2.8}),
the necessity of
(1)--(3) follows from \thmref{thm:strivial} and \corref{cor:symm}.
So assume that (1)--(3) hold.  If $\alpha\restr{\Om^\Sigma}$ is
implemented by $u$, then, to apply \thmref{thm:strivial}, we only
need to see that $\dtau^\w=u(\w,\cdot)$ is an (ordinary) Green
twisting map for the fiber $A_\w\rtimes_{\alpha^\w}G$ with domain
$\Sigma_\w$.  But this follows from
\cite[Lemma~6]{echros}\footnote{There are two ``Lemma~4''s in
\cite{echros}; presumably this reference will change.}.

Finally, the ``only if" direction of the last assertion is
immediate: if $\alpha$ is unitary, then $\acg\cong A\tensor
\cs(G)$.
On the other hand, if $\Prim(\acg)\cong\hat A\times\ghat G$, then
$\Sigma_\w^\perp=\set1$, and $\Om^\Sigma=\Om\times G$.
It is straightforward to check that $\alpha$ must be unitary if
the restriction of $\alpha$ to $\Om^\Sigma$ is unitary.
\end{pf}

It is interesting to consider a local version of \corref{cor:unitary}.
Recall that if $\alpha$ is locally unitary, then $\spec{\acg}$ is
a locally trivial principal $\ghat G$-bundle over $\hat A$
\cite{pr2}.
Conversely, if $\spec{\acg}$ is a locally trivial principal $\ghat
G$-bundle over $\hat A$, then we may apply \corref{cor:unitary}
to the local trivializations to conclude that $\alpha $ is
locally unitary.  In summary, we have the following corollary.
\begin{cor}
Suppose that $(A,\,G,\,\alpha)$ is a separable abelian dynamical
system with $\hat A$ Hausdorff.  Then $\spec{\acg} $ is a locally
trivial principal $\ghat G$-bundle over $\hat A$ if and only if
$\alpha$ is locally unitary.
\end{cor}

More generally, we can consider the case where $\Prim(\acg)$ is
any principal bundle.  Our result
generalizes \cite[Theorem~7.2]{90a}.
Here we do not assume that $A$ is type~I.
We omit the proof.
\begin{cor}\label{cor:trivialbundle}
Suppose that $(A,\,G,\,\alpha)$ is a separable abelian dynamical
system, and that $C$ is a closed subgroup of $G$.
Then, with respect to the dual action, $\Prim(\acg)$ is a trivial
$\ghat C$-bundle (that is, $\Prim(\acg)\cong_{\ghat
G}\Om\times\ghat C$) if
and only if the following three conditions are satisfied.
\begin{enumerate}
\item
The quasi-orbit space $\Om=\qob G[\Prim(A)]$ is Hausdorff.
\item
For all $\w\in\Om$, $C=\ttg(\alpha^\w)$.
\item
There is an (ordinary) Green twisting map $\tau$ for $(A,\,G,\,\alpha)$
with domain $C$.
\end{enumerate}
\end{cor}

Let us finally remark that our results here have
analogues in the case of abelian twisted covariant systems
$(A,\,G,\,\alpha,\,\tau)$ (i.e.\ $G/N_{\tau}$ is abelian).
These can be obtained by passing  to Morita equivalent ordinary
actions by applying \cite[Theorem 1]{ech5}.
To see that this works well, one only has to use the results about
Morita equivalent twisted actions as presented in \cite{ech5}.
We omit the straightforward details.

\def\mathcs{C^{\displaystyle *}} \def\cs{\ifmmode\mathcs\else$\mathcs$\fi}
\ifx\undefined\bysame
\newcommand{\bysame}{\leavevmode\hbox to3em{\hrulefill}\,}
\fi


\begin{thebibliography}{10}

\bibitem{ech4}
Siegfried Echterhoff, {\em Crossed products with continuous trace}, in
  preparation.

\bibitem{ech1}
\bysame, {\em The primitive ideal space of twisted covariant systems with
  continuously varying stabilizers}, Math. Ann. {\bf 292} (1992), 59--84.

\bibitem{ech5}
\bysame, {\em Morita equivalent actions and a new version of the
  {P}acker-{R}aeburn stabilization trick}, Jour. London Math. Soc. (1993), in
  press.

\bibitem{ech3}
\bysame, {\em On transformation group \cs-algebras with continuous trace},
  Trans. Amer. Math. Soc. (1993), in press.

\bibitem{ech2}
\bysame, {\em Regularizations of twisted covariant systems and crossed products
  with continuous trace}, J. Funct. Anal. (1993), in press.

\bibitem{echros}
Siegfried Echterhoff and Jonathan Rosenberg, {\em Fine structure of the
  {M}ackey machine for actions of abelian groups with constant {M}ackey
  obstuction}, preprint.

\bibitem{fell2}
J.~M.~G. Fell, {\em A {H}ausdorff topology on the closed subsets of a locally
  compact non-{H}ausdorff space}, Proc. Amer. Math. Soc. {\bf 13} (1962),
  472--476.

\bibitem{fell6}
\bysame, {\em Weak containment and induced representations of groups}, Canad.
  J. Math. {\bf 14} (1962), 237--268.

\bibitem{fell3}
\bysame, {\em Weak containment and induced representations of groups, {II}},
  Trans. Amer. Math. Soc. {\bf 110} (1964), 424--447.

\bibitem{fell4}
\bysame, {\em An extension of {M}ackey's method to {B}anach $*$-algebraic
  bundles}, Memoirs Amer. Math. Soc. {\bf 90} (1969), 1--168.

\bibitem{gl2}
James Glimm, {\em Families of induced representations}, Pacific J. Math. {\bf
  72} (1962), 885--911.

\bibitem{gl}
Elliot Gootman and Aldo Lazar, {\em Crossed products of type {I} {AF} algebras
  by abelian groups}, Israel J. Math. {\bf 56} (1986), 267--279.

\bibitem{goot2}
Elliot~C. Gootman, {\em Abelian group actions on type {I} {\cs}-algebras},
  Operator Algebras and their Connections with Topology and Ergodic Theory
  (Bu\c steni, Romania), Lecture Notes in Mathematics, vol. 1132,
  Springer-Verlag, 1983, pp.~152--169.

\bibitem{go}
Elliot~C. Gootman and Dorte Olesen, {\em Spectra of actions on type {I}
  {\cs}-algebras}, Math. Scand. {\bf 47} (1980), 329--349.

\bibitem{gr}
Elliot~C. Gootman and Jonathan Rosenburg, {\em The structure of crossed product
  {\cs}-algebras: A proof of the generalized {E}ffros-{H}ahn conjecture},
  Invent. Math. {\bf 52} (1979), 283--298.

\bibitem{green2}
Philip Green, {\em {\cs}-algebras of transformation groups with smooth orbit
  space}, Pacific J. Math. {\bf 72} (1977), 71--97.

\bibitem{green1}
\bysame, {\em The local structure of twisted covariance algebras}, Acta. Math.
  {\bf 140} (1978), 191--250.

\bibitem{horr}
Steven Hurder, Dorte Olesen, Iain Raeburn, and Jonathan Rosenberg, {\em The
  {C}onnes spectrum for actions of abelian groups on continuous-trace
  algebras}, Ergod. Th. \& Dynam. Sys. {\bf 6} (1986), 541--560.

\bibitem{lee2}
R.-Y. Lee, {\em On the \cs-algebras of operator fields}, Indiana Univ. Math. J.
  {\bf 25} (1976), 303--314.

\bibitem{op}
Dorte Olesen and Gert~K. Pedersen, {\em Partially inner {\cs}-dynamical
  systems}, J. Funct. Anal. {\bf 66} (1986), 263--281.

\bibitem{doir}
Dorte Olesen and Iain Raeburn, {\em Pointwise unitary automorphism groups}, J.
  Funct. Anal. {\bf 93} (1990), 278--309.

\bibitem{ped}
Gert~K. Pedersen, {\em {\cs}-algebras and their automorphism groups}, Academic
  Press, London, 1979.

\bibitem{pr2}
John Phillips and Iain Raeburn, {\em Crossed products by locally unitary
  automorphism groups and principal bundles}, J. Operator Theory {\bf 11}
  (1984), 215--241.

\bibitem{ra1}
Iain Raeburn, {\em On the {P}icard group of a continuous-trace {\cs}-algebra},
  Trans. Amer. Math. Soc. {\bf 263} (1981), 183--205.

\bibitem{rr}
Iain Raeburn and Jonathan Rosenberg, {\em Crossed products of continuous-trace
  {\cs}-algebras by smooth actions}, Trans. Amer. Math. Soc. {\bf 305} (1988),
  1--45.

\bibitem{rw2}
Iain Raeburn and Dana~P. Williams, {\em Crossed products by actions which are
  locally unitary on the stabilisers}, J. Funct. Anal. {\bf 81} (1988),
  385--431.

\bibitem{90a}
\bysame, {\em Moore cohomology, principal bundles, and actions of groups on
  {\cs}-algebras}, Indiana U. Math. J. {\bf 40} (1991), 707--740.

\bibitem{90b}
\bysame, {\em Topological invariants associated to the spectrum of crossed
  product {\cs}-algebras}, J. Funct. Anal. (1992), in press.

\bibitem{90c}
\bysame, {\em Dixmier-{D}ouady classes of dynamical systems and crossed
  products}, Canad. J. Math. (1993), in press.

\bibitem{95a}
\bysame, {\em Equivariant cohomology and a {G}ysin sequence for principal
  bundles}, preprint, March 1993.

\bibitem{schoch}
I.~Schochetman, {\em The dual topology of certain group extensions}, Adv. in
  Math. {\bf 35} (1980), 113--128.

\bibitem{dpw2}
Dana~P. Williams, {\em The topology on the primitive ideal space of
  transformation group {\cs}-algebras and {CCR} transformation group
  {\cs}-algebras}, Trans. Amer. Math. Soc. {\bf 266} (1981), 335--359.

\bibitem{dpw3}
\bysame, {\em Transformation group {\cs}-algebras with continuous trace}, J.
  Funct. Anal. {\bf 41} (1981), 40--76.

\bibitem{dpw4}
\bysame, {\em Transformation group {\cs}-algebras with {H}ausdorff spectrum},
  Illinois J. Math. {\bf 26} (1982), 317--321.

\bibitem{dpw5}
\bysame, {\em The structure of crossed products by smooth actions}, J. Austral.
  Math. Soc. (Series A) {\bf 47} (1989), 226--235.

\end{thebibliography}
\end{document}